\documentstyle[11pt]{article}

\newtheorem{eg}{{\bf Example}}
\newtheorem{lemma}{{\bf Lemma}}[section]
\newtheorem{theo}{{\bf Theorem}}
\newtheorem{prop}{{\bf Proposition}}[section]
\newtheorem{cor}[theo]{{\bf Corollary}}

\newtheorem{remark}{{\bf Remark}}

\font\bbb=msbm10 scaled\magstep1

\newcommand{\RR}{\mbox{\bbb R}}
\newcommand{\ZZ}{\mbox{\bbb Z}}

\def\scoll{\mbox{{\scriptsize
$\,\searrow$}$\!\!\stackrel{{}^{\rm s}}{}\,\,$}}

\def\coll{\mbox{\scriptsize $\,\,\searrow\,$}}

\newcommand{\bdN}{N^{\!\!\!^{^{\bullet}}}}

\newcommand{\si}{\sigma}

\begin{document}

\title{\bf Combinatorial triangulations of homology spheres}
\author{{\bf Bhaskar Bagchi}$^{\rm a}$, {\bf Basudeb Datta}$^{\rm b}$
}

\date{}

\maketitle


\noindent {\small $^{\rm a}$ Theoretical Statistics and
Mathematics Unit, Indian Statistical Institute,  Bangalore
560\,059, India.

\smallskip

\noindent $^{\rm b}$ Department of Mathematics, Indian Institute
of Science, Bangalore 560\,012,  India.}\footnotetext{{\em E-mail
addresses:} bbagchi@isibang.ac.in (B. Bagchi),
dattab@math.iisc.ernet.in (B. Datta).}

\begin{center}
\date{May 25, 2012}
\end{center}

\medskip

\hrule

\bigskip

 {\small

\noindent {\bf Abstract}

\bigskip

Let $M$ be an $n$-vertex combinatorial triangulation of a
$\ZZ_2$-homology $d$-sphere. In this paper we prove that if $n\leq
d+8$ then $M$ must be a combinatorial sphere. Further, if $n=d+9$
and $M$ is not a combinatorial sphere then $M$ can not admit any
proper bistellar move. Existence of a 12-vertex triangulation of
the lens space $L(3, 1)$ shows that the first result is sharp in
dimension three.

In the course of the proof we also show that any $\ZZ_2$-acyclic
simplicial complex on $\leq 7$ vertices is necessarily
collapsible. This result is best possible since there exist
$8$-vertex triangulations of the Dunce Hat which are not
collapsible. }

\bigskip

{\small


\smallskip


}

\bigskip

\hrule

\section{Introduction and results}

All the simplicial complexes considered in this paper are finite.
We say that a simplicial complex $K$ {\em triangulates} a
topological space $X$ (or $K$ is a {\em triangulation} of $X$) if
$X$ is homeomorphic to the geometric carrier $|K|$ of $K$.

The vertex-set of a simplicial complex $K$ is denoted by $V(K)$.
If $K$, $L$ are two simplicial complexes, then a {\em simplicial
isomorphism} from $K$ to $L$ is a bijection $\pi : V(K) \to V(L)$
such that for $\si\subseteq V(K)$, $\si$ is a face of $K$ if and
only if $\pi (\si)$ is a face of $L$. The complexes $K$, $L$ are
called (simplicially) {\em isomorphic} when such an isomorphism
exists. We identify two simplicial complexes if they are
isomorphic.

A simplicial complex $K$ is called {\em pure} if all the maximal
faces of $K$ have the same dimension.  A maximal face in a pure
simplicial complex is also called a {\em facet}.

If $\si$ is a face of a simplicial complex  $K$ then the {\em
link} of $\si$ in $K$, denoted by ${\rm Lk}_K(\sigma)$ (or simply
by ${\rm Lk}(\sigma)$), is by definition the simplicial complex
whose faces are the faces $\tau$ of $K$ such that $\tau$ is
disjoint from $\si$ and $\si\cup\tau$ is a face of $K$.

A subcomplex $L$ of a simplicial complex $K$ is called an {\em
induced} (or {\em full}\,) subcomplex of $K$ if $\sigma\in K$ and
$\sigma\subseteq V(L)$ imply $\sigma \in L$. The induced
subcomplex of $K$ on the vertex set $U$ is denoted by $K[U]$.

For a commutative ring $R$, a simplicial complex $K$ is called
$R$-acyclic if $|K|$ is $R$-acyclic, i.e., $\widetilde{H}_q(|K|,
R) = 0$ for all $q\geq 0$ (where $\widetilde{H}^q(|K|, R)$ denotes
the reduced homology).

By a {\em subdivision} of a simplicial complex  $K$ we mean a
simplicial complex $K^{\,\prime}$ together with a homeomorphism
from $|K^{\,\prime}|$ onto $|K|$ which is facewise linear. Two
simplicial complexes $K$ and $L$ are called {\em combinatorially
equivalent} (denoted by {\em $K \approx L$}) if they have
isomorphic subdivisions. So, $K\approx L$ if and only if $|K|$ and
$|L|$ are piecewise-linear (pl) homeomorphic (see \cite{rs}).

For a set $U$ with $d + 1$ elements, let $K$ be the simplicial
complex whose faces are all the non-empty subsets of $U$. Then $K$
triangulates the $d$-dimensional closed unit ball. This complex is
called the {\em standard $d$-ball} and is denoted by
$\Delta^{d}_{d + 1}(U)$ or simply by $\Delta^{d}_{d + 1}$. A
polyhedron is called a {\em pl $d$-ball} if it is pl homeomorphic
to $|\Delta^d_{d+1}|$. A simplicial complex $X$ is called a {\em
combinatorial $d$-ball} if it is combinatorially equivalent to
$\Delta^{d}_{d+ 1}$. So, $X$ is a combinatorial $d$-ball if and
only if $|X|$ is a pl $d$-ball.

For a set $V$ with $d + 2$ elements, let $S$ be the simplicial
complex whose faces are all the non-empty proper subsets of $V$.
Then $S$ triangulates the $d$-sphere. This complex is called the
{\em standard $d$-sphere} and is denoted by $S^{\,d}_{d + 2}(V)$
or simply by $S^{\,d}_{d + 2}$. A polyhedron is called a {\em pl
$d$-sphere} if it is pl homeomorphic to $|S^{\,d}_{d + 2}|$. A
simplicial complex  $X$ is called a {\em combinatorial $d$-sphere}
if it is combinatorially equivalent to $S^{\,d}_{d + 2}$.  So, $X$
is a combinatorial $d$-sphere if and only if $|X|$ is a pl
$d$-sphere.

A simplicial complex $K$ is called a {\em combinatorial
$d$-manifold} if the link of each vertex is a combinatorial
$(d-1)$-sphere. A simplicial complex  $K$ is a combinatorial
$d$-manifold if and only if $|K|$ is a closed pl $d$-manifold (see
\cite{rs}).

If a triangulation $K$ of a space $X$ is a combinatorial manifold
then $K$ is called a {\em combinatorial triangulation} of $X$. If
$K$ is a triangulation of a 3-manifold then the link of a vertex
is a triangulation of the 2-sphere and all triangulations of the
2-sphere are combinatorial 2-spheres. So, any triangulation of a
3-manifold is a combinatorial triangulation.

Let $\tau\subset\sigma$ be two faces of a simplicial complex $K$.
We say that $\tau$ is a {\em free face} of $\sigma$ if $\sigma$ is
the only face of $K$ which properly contains $\tau$. (It follows
that $\dim(\sigma)-\dim(\tau)=1$ and $\sigma$ is a maximal simplex
in $K$.) If $\tau$ is a free face of $\sigma$ then $K^{\,\prime}
:= K \setminus \{\tau, \sigma\}$ is a simplicial complex. We say
that there is an {\em elementary collapse} of $K$ to
$K^{\,\prime}$. We say $K$ {\em collapses} to $L$ and write
$K\scoll L$ if there exists a sequence $K=K_0$, $K_1, \dots$,
$K_n=L$ of simplicial complexes such that there is an elementary
collapse of $K_{i-1}$ to $K_{i}$ for $1\leq i\leq n$ (see
\cite{b}). If $L$ consists of a 0-simplex (a point) we say that
$K$ is {\em collapsible} and write $K\scoll 0$. Clearly, if
$K\scoll L$ then $|K|\coll |L|$ as polyhedra and hence $|K|$ and
$|L|$ have the same homotopy type (see \cite{rs}). So, if a
simplicial complex $K$ is collapsible then $|K|$ is contractible
and hence, in particular, $K$ is $\ZZ_2$-acyclic. Here we prove\,:

\begin{theo}$\!\!${\bf .} \label{t1}
If a $\ZZ_2$-acyclic simplicial complex has $\leq 7$ vertices then
it is collapsible.
\end{theo}

As an application of Theorem \ref{t1}, we prove our main result -
a recognition theorem for combinatorial spheres\,:

\begin{theo}$\!\!${\bf .} \label{t2}
Let $M$ be an $n$-vertex combinatorial triangulation of a
$\ZZ_2$-homology $d$-sphere. Suppose $M$ has an $m$-vertex
combinatorial $d$-ball as an induced subcomplex, where $n \leq m +
7$. Then $M$ is a combinatorial sphere.
\end{theo}

In consequence we get the following.

\begin{cor}$\!\!${\bf .} \label{t3}
Let $M$ be an $n$-vertex combinatorial $d$-manifold. If $|M|$ is a
$\ZZ_2$-homology sphere and $n \leq d+8$ then $M$ is a
combinatorial sphere.
\end{cor}

\begin{cor}$\!\!${\bf .} \label{t4}
Let $M$ be a $(d + 9)$-vertex combinatorial triangulation of a
$\ZZ_2$-homology $d$- sphere. If $M$ is not a combinatorial sphere
then $M$ can not admit any bistellar $i$-move for $i < d$.
\end{cor}

Since by the universal coefficient theorem any integral homology
sphere is a $\ZZ_2$-homology sphere, Theorem \ref{t2}, Corollary
\ref{t3} and Corollary \ref{t4} remain true if we replace
$\ZZ_2$-homology by integral homology in the hypothesis. In
particular, we have\,:

\begin{cor}$\!\!${\bf .} \label{t5}
Let $M$ be an $n$-vertex combinatorial triangulation of an
integral homology $d$-sphere.
\begin{enumerate}
   \vspace*{-1mm} \item[{$(a)$}] If $n \leq d+8$ then $M$ is a
combinatorial sphere.
   \vspace*{-1.25mm} \item[{$(b)$}] If $n = d + 9$ and $M$ is not
a combinatorial sphere then $M$ can not admit any bistellar
$i$-move for $i < d$.
\end{enumerate}
\end{cor}

\begin{remark}$\!\!${\bf .}
{\rm Corollary \ref{t3} is clearly trivial for $d\leq 2$. In
\cite{bk}, Brehm and K\"{u}hnel proved that any $n$-vertex
combinatorial $d$-manifold is a combinatorial $d$-sphere if $n < 3
\lceil d/2 \rceil + 3$ and it is either a combinatorial $d$-sphere
or a cohomology projective plane if $n = 3d/2 + 3$. So, Corollary
\ref{t3} has new content only for $3 \leq d \leq 8$. }
\end{remark}

\begin{remark}$\!\!${\bf .}
{\rm Another result in \cite{bk} says that any $n$-vertex
combinatorial $d$-manifold is simply connected for $n\leq 2d + 2$.
Since a simply connected integral homology sphere is a sphere for
$d\not = 3$, and since for $d\not = 4$ all combinatorial
triangulations of $d$-spheres are combinatorial spheres, this
result implies that all combinatorial triangulations of integral
homology $d$-spheres ($d\not = 3, 4$) with $\leq 2d+2$ vertices
are combinatorial spheres. This is stronger than Corollary
\ref{t5} $(a)$ for $d\geq 6$. Thus Corollary \ref{t5} $(a)$ has
new content only for $d = 3, 4, 5$. }
\end{remark}

\begin{remark}$\!\!${\bf .}
{\rm In \cite[p. 35]{l}, Lutz presented a 12-vertex combinatorial
triangulation of the lens space $L(3, 1)$. (It is mentioned in
\cite[p. 79]{k} that Brehm obtained a 12-vertex combinatorial
triangulation of $L(3, 1)$ earlier.) Since $L(3, 1)$ is a
$\ZZ_2$-homology $3$-sphere ($H_1(L(3,1), \ZZ)=\ZZ_3$,
$H_2(L(3,1), \ZZ)= 0$), Corollary \ref{t3} is sharp for $d = 3$.

It follows from Corollary \ref{t3} that 12 is the least number of
vertices required to triangulate $L(3, 1)$. It follows from
Corollary \ref{t4} that a 12-vertex combinatorial triangulation of
$L(3, 1)$ can not admit any bistellar $i$-move for $0\leq i\leq
2$.  }
\end{remark}

\begin{remark}$\!\!${\bf .}
{\rm Recall that the Dunce Hat is the topological space obtained
from the solid triangle $abc$ by identifying the oriented edges
$\vec{ab}$, $\vec{bc}$ and $\vec{ac}$. The following is a
triangulation of the Dunce Hat using 8 vertices.}

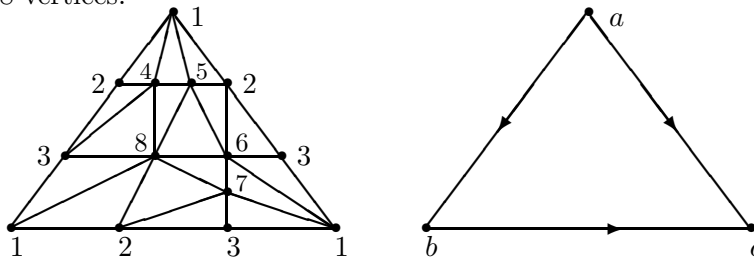
\setlength{\unitlength}{2.4mm}
\begin{picture}(47,14)(-13,0)

\put(1.77,2){$_{\bullet}$} \put(7.77,2){$_{\bullet}$}
\put(13.77,2){$_{\bullet}$} \put(19.77,2){$_{\bullet}$}

\put(13.77,4){$_{\bullet}$} \put(10.77,14){$_{\bullet}$}

\put(4.77,6){$_{\bullet}$} \put(9.77,6){$_{\bullet}$}
\put(13.77,6){$_{\bullet}$} \put(16.77,6){$_{\bullet}$}

\put(7.77,10){$_{\bullet}$} \put(9.77,10){$_{\bullet}$}
\put(11.77,10){$_{\bullet}$} \put(13.77,10){$_{\bullet}$}

\put(2,0.5){$1$} \put(8,0.5){$2$} \put(14,0,5){$3$}
\put(20,0.5){$1$} \put(3.5,5.5){$3$} \put(6.5,9.5){$2$}
\put(17.9,5.5){$3$} \put(14.9,9.5){$2$} \put(12,13.2){$1$}

\put(8.9,6.3){\small 8} \put(14.5,6.3){\small 6}
\put(14.5,4.1){\footnotesize 7} \put(9.2,10.3){\footnotesize 4}
\put(12.3,10.3){\footnotesize 5}

 \thicklines

\put(2,2){\line(1,0){18}} \put(2,2){\line(3,4){9}}
\put(20,2){\line(-3,4){9}} \put(10,10){\line(1,4){1}}
\put(14,2){\line(0,1){8}} \put(5,6){\line(1,0){12}}
\put(8,10){\line(1,0){6}}

\put(2,2){\line(2,1){8}} \put(10,6){\line(0,1){4}}
\put(8,2){\line(1,2){2}} \put(8,2){\line(3,1){6}}
\put(20,2){\line(-3,1){6}} \put(10,6){\line(2,-1){4}}
\put(14,6){\line(3,-2){6}} \put(5,6){\line(5,4){5}}
\put(10,6){\line(1,2){2}} \put(14,6){\line(-1,2){2}}

\put(12,10){\line(-1,4){1}}


\put(24.77,2){$_{\bullet}$} \put(42.77,2){$_{\bullet}$}

\put(33.77,14){$_{\bullet}$}

\put(25,0.5){$b$} \put(43,0.5){$c$} \put(35.2,13.2){$a$}

 \thicklines

\put(25,2){\line(1,0){18}} \put(25,2){\line(3,4){9}}
\put(43,2){\line(-3,4){9}}

\put(31,10){\vector(-3,-4){2}} \put(37,10){\vector(3,-4){2}}
\put(34,2){\vector(1,0){2}}

\end{picture}

\noindent {\rm Since this example is contractible but not
collapsible, it follows that the bound 7 in Theorem \ref{t1} is
best possible.}
\end{remark}

\begin{remark}$\!\!${\bf .}
{\rm Let $H^{\,3}$ be the non-orientable 3-manifold obtained from
$S^{\,2}\times [0, 1]$ by identifying $(x, 0)$ with $(-x, 1)$. It
follows from works of Walkup \cite[Theorems 3, 4]{w} that} if $K$
is a combinatorial $3$-manifold and $|K|$ is not homeomorphic to
$S^{\,3}$, $S^{\,2} \times S^{\,1}$ or $H^{\,3}$ then $f_1(K) \geq
4f_0(K) + 8$ and hence $f_0(K)\geq 11$. {\rm Thus if $M$ ($\neq
S^{\,3}$) is a $\ZZ_2$-homology 3-sphere then at least 11 vertices
are needed for any combinatorial triangulation of $M$. Now,
Corollary \ref{t3} implies that at least 12 vertices are needed.
In \cite{bl}, Bj\"{o}rner and Lutz have presented a 16-vertex
combinatorial triangulation of the Poincar\'{e} homology 3-sphere.

In \cite{bd4}, we have shown that all combinatorial triangulations
of $S^{\,4}$ with at most 10 vertices are combinatorial
$4$-spheres. Now, Corollary \ref{t3} implies that all
combinatorial triangulations of $S^{\,4}$ with at most 12 vertices
are combinatorial spheres. So, any combinatorial triangulation (if
it exists) of $S^{\,4}$ which is not a combinatorial sphere
requires at least 13 vertices.}
\end{remark}

\begin{remark}$\!\!${\bf .}
{\rm The conclusion in Corollary \ref{t4} (namely, that certain
combinatorial manifolds do not admit any proper bistellar move)
appears to be a strong structural restriction. We owe to F. H.
Lutz the information that the smallest known combinatorial sphere
(other than a standard sphere) not admitting any proper bistellar
move is a 16-vertex 3-sphere. }
\end{remark}

\section{Preliminaries and Definitions.}

For a simplicial complex $K$, the maximum $k$ such that $K$ has a
$k$-face is called the {\em dimension} of $K$. An 1-dimensional
simplicial complex is called a {\em graph}. A simplicial complex
$K$ is called {\em connected} if $|K|$ is connected.

For $i = 1, 2, 3$, the $i$-faces of a simplicial complex are also
called the {\em edges}, {\em triangles} and {\em tetrahedra} of
the complex, respectively. For a face $\sigma$ in a simplicial
complex $K$, the number of vertices in ${\rm Lk}_K(\sigma)$ is
called the {\em degree} of $\sigma$ in $K$ and is denoted by
$\deg_K(\sigma)$.

If the number of $i$-simplices of a $d$-dimensional simplicial
complex  $K$ is $f_i(K)$, then the vector $f= (f_0, \dots, f_d)$
is called the {\em $f$-vector} of $K$ and the number $\chi(K) :=
\sum_{i= 0}^{d}(- 1)^i f_i(K)$ is called the {\em Euler
characteristic}  of $K$. If $f_{k-1}= {f_0 \choose k}$ then $K$ is
called {\em $k$-neighbourly}.

For two simplicial complexes $K$, $L$ with disjoint vertex sets,
the {\em join} $K \ast L$ is the simplicial complex $K\cup L\cup
\{\sigma \cup \tau ~ : ~ \sigma \in K, \tau\in L\}$.

If $K$ is a $d$-dimensional simplicial complex then define the
{\em pure part} of $K$ as the simplicial complex whose simplices
are the sub-simplices of the $d$-simplices of $K$.

A $d$-dimensional pure simplicial complex $K$ is called a {\em
weak pseudomanifold} if each $(d-1)$-face is contained in exactly
two facets of $K$. A $d$-dimensional weak pseudomanifold $K$ is
called a {\em pseudomanifold} if for any pair $\tau$, $\sigma$ of
facets, there exists a sequence $\tau = \tau_0, \dots,
\tau_n=\sigma$ of facets of $K$, such that $\tau_{i - 1} \cap
\tau_{i}$ is a $(d-1)$-simplex of $K$ for $1\leq i\leq n$. In
other words, a weak pseudomanifold is a pseudomanifold if and only
if it does not have any weak pseudomanifold of the same dimension
as a proper subcomplex. Clearly, any connected combinatorial
manifold is a pseudomanifold.

For $n\geq 3$, the $n$-vertex combinatorial 1-sphere ({\em
$n$-cycle}) is the unique $n$-vertex 1-dimensional  pseudomanifold
and is denoted by $S^{\,1}_n$.

A $d$-dimensional pure simplicial complex $K$ is called a {\em
weak pseudomanifold with boundary} if each $(d-1)$-face is
contained in 1 or 2 facets of $K$ and there exists a $(d-1)$-face
of degree 1. The boundary $\partial K$ of $K$ is by definition the
pure simplicial complex whose facets are the degree one
$(d-1)$-faces of $K$.

A simplicial complex $K$ is called a {\em combinatorial
$d$-manifold with boundary} if the link of each vertex is either a
combinatorial $(d-1)$-sphere or a combinatorial $(d-1)$-ball and
there exists a vertex whose link is a combinatorial $(d-1)$-ball.
A simplicial complex $K$ is a combinatorial $d$-manifold with
boundary if and only if $|K|$ is a compact pl $d$-manifold with
non-empty boundary. Clearly, if $K$ is a combinatorial
$d$-manifold with boundary then $\partial K\neq\emptyset$ and
${\rm Lk}_{\partial K}(v) = \partial({\rm Lk}_K(v))$, for $v\in
V(\partial K)$. Therefore, $\partial K$ is a combinatorial
$(d-1)$-manifold. Clearly, if $K$ is a combinatorial $d$-ball
($d>0$) then $K$ is a combinatorial $d$-manifold with boundary and
$\partial K$ is a combinatorial $(d-1)$-sphere.

\begin{eg}$\!\!${\bf .} \label{e1}
{\rm Some weak pseudomanifolds on 6 or 7 vertices.}


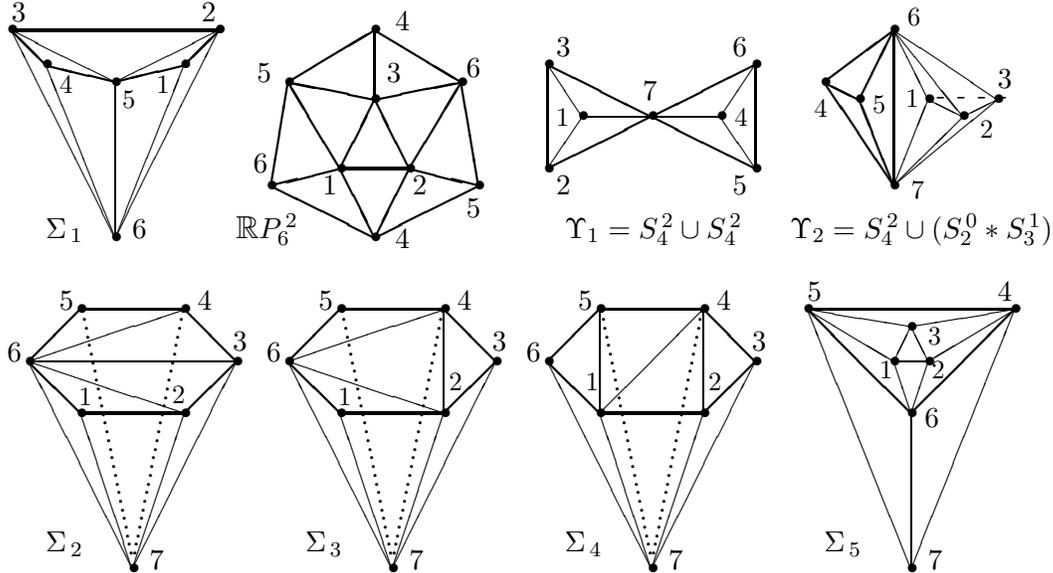
\setlength{\unitlength}{2.3mm}
\begin{picture}(60,15)(0,-1)


\put(6.77,9){$_{\bullet}$} \put(2.77,10){$_{\bullet}$}
\put(0.77,12){$_{\bullet}$} \put(12.77,12){$_{\bullet}$}
\put(10.77,10){$_{\bullet}$} \put(6.77,0){$_{\bullet}$}

\put(9.4,8.4){\small 1} \put(12,12.5){$2$} \put(1,12.5){$3$}
\put(3.8,8.4){\small 4} \put(7.4,7.5){\small 5} \put(8,0){$6$}

\thicklines

\put(7,9){\line(-4,1){4}} \put(3,10){\line(-1,1){2}}
\put(1,12){\line(1,0){12}} \put(13,12){\line(-1,-1){2}}
\put(11,10){\line(-4,-1){4}}

\thinlines

\put(7,9){\line(-2,1){6}} \put(7,9){\line(2,1){6}}
\put(7,0){\line(-1,2){6}} \put(7,0){\line(-2,5){4}}
\put(7,0){\line(0,1){9}} \put(7,0){\line(2,5){4}}
\put(7,0){\line(1,2){6}}

\put(3,0){\mbox{$\Sigma_{\,1}$}}


\put(21.77,0){$_{\bullet}$}
\put(15.77,3){$_{\bullet}$}    
\put(16.77,9){$_{\bullet}$} \put(21.77,12){$_{\bullet}$}
\put(26.77,9){$_{\bullet}$} \put(27.77,3){$_{\bullet}$}
\put(19.77,4){$_{\bullet}$} \put(23.77,4){$_{\bullet}$}
\put(21.77,8){$_{\bullet}$}

\put(23.2,-0.6){$4$} \put(15,3.8){$6$} \put(15.2,9){$5$}
\put(23.2,12){$4$} \put(27.3,9.3){$6$} \put(27.3,1.2){$5$}
\put(19,2.5){$1$} \put(24.2,2.5){$2$} \put(22.7,9){$3$}

\thicklines

\put(22,0){\line(-2,1){6}} \put(16,3){\line(1,6){1}}
\put(17,9){\line(5,3){5}} \put(22,12){\line(5,-3){5}}
\put(27,9){\line(1,-6){1}} \put(28,3){\line(-2,-1){6}}

\put(20,4){\line(-4,-1){4}} \put(20,4){\line(-3,5){3}}
\put(20,4){\line(1,2){2}} \put(20,4){\line(1,0){4}}
\put(20,4){\line(1,-2){2}}

\put(24,4){\line(-1,-2){2}} \put(24,4){\line(-1,2){2}}
\put(24,4){\line(3,5){3}} \put(24,4){\line(4,-1){4}}

\put(22,8){\line(-5,1){5}} \put(22,8){\line(0,1){4}}
\put(22,8){\line(5,1){5}}

\put(14,0){\mbox{$\RR P^{\,2}_6$}}



\put(31.77,4){$_{\bullet}$} \put(43.77,4){$_{\bullet}$}
\put(31.77,10){$_{\bullet}$} \put(43.77,10){$_{\bullet}$}
\put(33.77,7){$_{\bullet}$} \put(41.77,7){$_{\bullet}$}
\put(37.77,7){$_{\bullet}$}

\put(32.5,6.5){\small 1} \put(32.5,2.5){$2$} \put(32.5,10.5){$3$}
\put(42.8,6.5){\small 4} \put(42.7,2.5){$5$} \put(42.7,10.5){$6$}
\put(37.5,8){$7$}

\thicklines

\put(32,4){\line(2,1){6}} \put(32,4){\line(0,1){6}}
\put(44,4){\line(-2,1){6}} \put(44,4){\line(0,1){6}}
\put(32,10){\line(2,-1){6}} \put(44,10){\line(-2,-1){6}}

\thinlines

\put(32,4){\line(2,3){2}} \put(32,10){\line(2,-3){2}}
\put(44,4){\line(-2,3){2}} \put(44,10){\line(-2,-3){2}}
\put(34,7){\line(1,0){8}}

\put(33,0){\mbox{$\Upsilon_1=S^{\,2}_4\cup S^{\,2}_4$}}


\put(51.77,3){$_{\bullet}$} \put(55.77,7){$_{\bullet}$}
\put(49.77,8){$_{\bullet}$} \put(53.77,8){$_{\bullet}$}
\put(57.77,8){$_{\bullet}$} \put(47.77,9){$_{\bullet}$}
\put(51.77,12){$_{\bullet}$}

\put(54.35,7.78){- - - -} 

\put(52.7,7.6){\small 1} \put(57,5.8){\small 2} \put(58,8.7){$3$}
\put(47.4,7){\small $4$} \put(50.8,7.2){\small 5}
\put(52.8,12.3){\small $6$} \put(53,2.3){\small $7$}

\thicklines

\put(48,9){\line(4,3){4}} \put(50,8){\line(1,2){2}}
\put(50,8){\line(-2,1){2}} \put(52,3){\line(-2,3){4}}
\put(52,3){\line(-2,5){2}} \put(52,3){\line(0,1){9}}

\thinlines

\put(56,7){\line(2,1){2}} \put(56,7){\line(-2,1){2}}

\put(52,3){\line(2,5){2}} \put(52,3){\line(1,1){4}}
\put(52,3){\line(6,5){6}} \put(52,12){\line(1,-2){2}}
\put(52,12){\line(4,-5){4}} \put(52,12){\line(3,-2){6}}

\put(46,0){\mbox{$\Upsilon_2 = S^{\,2}_4\cup (S^{\,0}_2\ast
S^{\,1}_3)$}}


\end{picture}

\setlength{\unitlength}{2.3mm}
\begin{picture}(60,20)(-2,-2)


\put(2.77,9){$_{\bullet}$} \put(-0.23,12){$_{\bullet}$}
\put(2.77,15){$_{\bullet}$} \put(8.77,15){$_{\bullet}$}
\put(11.77,12){$_{\bullet}$} \put(8.77,9){$_{\bullet}$}
\put(5.77,0){$_{\bullet}$}

\put(2.9,9.5){\small $1$} \put(8.3,9.6){\small 2}
\put(11.8,12.6){$3$} \put(9.77,15){$4$} \put(1.77,15){$5$}
\put(-1.23,12){$6$}
  \put(7,0){$7$}

\thicklines

\put(3,9){\line(-1,1){3}} \put(0,12){\line(1,1){3}}
\put(3,15){\line(1,0){6}} \put(9,15){\line(1,-1){3}}
\put(12,12){\line(-1,-1){3}} \put(9,9){\line(-1,0){6}}

\thinlines

\put(0,12){\line(1,0){12}} \put(0,12){\line(3,1){9}}
\put(0,12){\line(3,-1){9}}

\put(6,0){\line(-1,2){6}} \put(6,0){\line(1,2){6}}
\put(6,0){\line(-1,3){3}} \put(6,0){\line(1,3){3}}

\put(5.7,0.5){\mbox{$\cdot$}} \put(5.6,1){\mbox{$\cdot$}}
\put(5.5,1.5){\mbox{$\cdot$}} \put(5.4,2){\mbox{$\cdot$}}
\put(5.3,2.5){\mbox{$\cdot$}} \put(5.2,3){\mbox{$\cdot$}}
\put(5.1,3.5){\mbox{$\cdot$}} \put(5,4){\mbox{$\cdot$}}
\put(4.9,4.5){\mbox{$\cdot$}} \put(4.8,5){\mbox{$\cdot$}}
\put(4.7,5.5){\mbox{$\cdot$}} \put(4.6,6){\mbox{$\cdot$}}
\put(4.5,6.5){\mbox{$\cdot$}} \put(4.4,7){\mbox{$\cdot$}}
\put(4.3,7.5){\mbox{$\cdot$}} \put(4.2,8){\mbox{$\cdot$}}
\put(4.1,8.5){\mbox{$\cdot$}} \put(4,9){\mbox{$\cdot$}}
\put(3.9,9.5){\mbox{$\cdot$}} \put(3.8,10){\mbox{$\cdot$}}
\put(3.7,10.5){\mbox{$\cdot$}} \put(3.6,11){\mbox{$\cdot$}}
\put(3.5,11.5){\mbox{$\cdot$}} \put(3.4,12){\mbox{$\cdot$}}
\put(3.3,12.5){\mbox{$\cdot$}} \put(3.2,13){\mbox{$\cdot$}}
\put(3.1,13.5){\mbox{$\cdot$}} \put(3,14){\mbox{$\cdot$}}
\put(2.9,14.5){\mbox{$\cdot$}}

\put(5.9,0.5){\mbox{$\cdot$}} \put(6,1){\mbox{$\cdot$}}
\put(6.1,1.5){\mbox{$\cdot$}} \put(6.2,2){\mbox{$\cdot$}}
\put(6.3,2.5){\mbox{$\cdot$}} \put(6.4,3){\mbox{$\cdot$}}
\put(6.5,3.5){\mbox{$\cdot$}} \put(6.6,4){\mbox{$\cdot$}}
\put(6.7,4.5){\mbox{$\cdot$}} \put(6.8,5){\mbox{$\cdot$}}
\put(6.9,5.5){\mbox{$\cdot$}} \put(7,6){\mbox{$\cdot$}}
\put(7.1,6.5){\mbox{$\cdot$}} \put(7.2,7){\mbox{$\cdot$}}
\put(7.3,7.5){\mbox{$\cdot$}} \put(7.4,8){\mbox{$\cdot$}}
\put(7.5,8.5){\mbox{$\cdot$}} \put(7.6,9){\mbox{$\cdot$}}
\put(7.7,9.5){\mbox{$\cdot$}} \put(7.8,10){\mbox{$\cdot$}}
\put(7.9,10.5){\mbox{$\cdot$}} \put(8,11){\mbox{$\cdot$}}
\put(8.1,11.5){\mbox{$\cdot$}} \put(8.2,12){\mbox{$\cdot$}}
\put(8.3,12.5){\mbox{$\cdot$}} \put(8.4,13){\mbox{$\cdot$}}
\put(8.5,13.5){\mbox{$\cdot$}} \put(8.6,14){\mbox{$\cdot$}}
\put(8.7,14.5){\mbox{$\cdot$}}

\put(1,1){\mbox{$\Sigma_{\,2}$}}


\put(17.77,9){$_{\bullet}$} \put(14.77,12){$_{\bullet}$}
\put(17.77,15){$_{\bullet}$} \put(23.77,15){$_{\bullet}$}
\put(26.77,12){$_{\bullet}$} \put(23.77,9){$_{\bullet}$}
\put(20.77,0){$_{\bullet}$}

\put(17.9,9.5){\small 1} \put(24.3,10.3){\small 2}
\put(26.8,12.6){$3$} \put(24.77,15){$4$} \put(16.77,15){$5$}
\put(13.77,12){$6$}
  \put(22,0){$7$}

\thicklines

\put(18,9){\line(-1,1){3}} \put(15,12){\line(1,1){3}}
\put(18,15){\line(1,0){6}} \put(24,15){\line(1,-1){3}}
\put(27,12){\line(-1,-1){3}} \put(24,9){\line(-1,0){6}}

\thinlines

\put(24,9){\line(0,1){6}} \put(15,12){\line(3,1){9}}
\put(15,12){\line(3,-1){9}}

\put(21,0){\line(-1,2){6}} \put(21,0){\line(1,2){6}}
\put(21,0){\line(-1,3){3}} \put(21,0){\line(1,3){3}}

\put(20.7,0.5){\mbox{$\cdot$}} \put(20.6,1){\mbox{$\cdot$}}
\put(20.5,1.5){\mbox{$\cdot$}} \put(20.4,2){\mbox{$\cdot$}}
\put(20.3,2.5){\mbox{$\cdot$}} \put(20.2,3){\mbox{$\cdot$}}
\put(20.1,3.5){\mbox{$\cdot$}} \put(20,4){\mbox{$\cdot$}}
\put(19.9,4.5){\mbox{$\cdot$}} \put(19.8,5){\mbox{$\cdot$}}
\put(19.7,5.5){\mbox{$\cdot$}} \put(19.6,6){\mbox{$\cdot$}}
\put(19.5,6.5){\mbox{$\cdot$}} \put(19.4,7){\mbox{$\cdot$}}
\put(19.3,7.5){\mbox{$\cdot$}} \put(19.2,8){\mbox{$\cdot$}}
\put(19.1,8.5){\mbox{$\cdot$}} \put(19,9){\mbox{$\cdot$}}
\put(18.9,9.5){\mbox{$\cdot$}} \put(18.8,10){\mbox{$\cdot$}}
\put(18.7,10.5){\mbox{$\cdot$}} \put(18.6,11){\mbox{$\cdot$}}
\put(18.5,11.5){\mbox{$\cdot$}} \put(18.4,12){\mbox{$\cdot$}}
\put(18.3,12.5){\mbox{$\cdot$}} \put(18.2,13){\mbox{$\cdot$}}
\put(18.1,13.5){\mbox{$\cdot$}} \put(18,14){\mbox{$\cdot$}}
\put(17.9,14.5){\mbox{$\cdot$}}

\put(20.9,0.5){\mbox{$\cdot$}} \put(21,1){\mbox{$\cdot$}}
\put(21.1,1.5){\mbox{$\cdot$}} \put(21.2,2){\mbox{$\cdot$}}
\put(21.3,2.5){\mbox{$\cdot$}} \put(21.4,3){\mbox{$\cdot$}}
\put(21.5,3.5){\mbox{$\cdot$}} \put(21.6,4){\mbox{$\cdot$}}
\put(21.7,4.5){\mbox{$\cdot$}} \put(21.8,5){\mbox{$\cdot$}}
\put(21.9,5.5){\mbox{$\cdot$}} \put(22,6){\mbox{$\cdot$}}
\put(22.1,6.5){\mbox{$\cdot$}} \put(22.2,7){\mbox{$\cdot$}}
\put(22.3,7.5){\mbox{$\cdot$}} \put(22.4,8){\mbox{$\cdot$}}
\put(22.5,8.5){\mbox{$\cdot$}} \put(22.6,9){\mbox{$\cdot$}}
\put(22.7,9.5){\mbox{$\cdot$}} \put(22.8,10){\mbox{$\cdot$}}
\put(22.9,10.5){\mbox{$\cdot$}} \put(23,11){\mbox{$\cdot$}}
\put(23.1,11.5){\mbox{$\cdot$}} \put(23.2,12){\mbox{$\cdot$}}
\put(23.3,12.5){\mbox{$\cdot$}} \put(23.4,13){\mbox{$\cdot$}}
\put(23.5,13.5){\mbox{$\cdot$}} \put(23.6,14){\mbox{$\cdot$}}
\put(23.7,14.5){\mbox{$\cdot$}}

\put(16,1){\mbox{$\Sigma_{\,3}$}}


\put(32.77,9){$_{\bullet}$} \put(29.77,12){$_{\bullet}$}
\put(32.77,15){$_{\bullet}$} \put(38.77,15){$_{\bullet}$}
\put(41.77,12){$_{\bullet}$} \put(38.77,9){$_{\bullet}$}
\put(35.77,0){$_{\bullet}$}

\put(32.2,10.3){\small 1} \put(39.3,10.3){\small 2}
\put(41.8,12.6){$3$} \put(39.77,15){$4$} \put(31.77,15){$5$}
\put(28.77,12){$6$}
  \put(37,0){$7$}

\thicklines

\put(33,9){\line(-1,1){3}} \put(30,12){\line(1,1){3}}
\put(33,15){\line(1,0){6}} \put(39,15){\line(1,-1){3}}
\put(42,12){\line(-1,-1){3}} \put(39,9){\line(-1,0){6}}

\thinlines

\put(33,9){\line(0,1){6}} \put(33,9){\line(1,1){6}}
\put(39,9){\line(0,1){6}}

\put(36,0){\line(-1,2){6}} \put(36,0){\line(1,2){6}}
\put(36,0){\line(-1,3){3}} \put(36,0){\line(1,3){3}}

\put(35.7,0.5){\mbox{$\cdot$}} \put(35.6,1){\mbox{$\cdot$}}
\put(35.5,1.5){\mbox{$\cdot$}} \put(35.4,2){\mbox{$\cdot$}}
\put(35.3,2.5){\mbox{$\cdot$}} \put(35.2,3){\mbox{$\cdot$}}
\put(35.1,3.5){\mbox{$\cdot$}} \put(35,4){\mbox{$\cdot$}}
\put(34.9,4.5){\mbox{$\cdot$}} \put(34.8,5){\mbox{$\cdot$}}
\put(34.7,5.5){\mbox{$\cdot$}} \put(34.6,6){\mbox{$\cdot$}}
\put(34.5,6.5){\mbox{$\cdot$}} \put(34.4,7){\mbox{$\cdot$}}
\put(34.3,7.5){\mbox{$\cdot$}} \put(34.2,8){\mbox{$\cdot$}}
\put(34.1,8.5){\mbox{$\cdot$}} \put(34,9){\mbox{$\cdot$}}
\put(33.9,9.5){\mbox{$\cdot$}} \put(33.8,10){\mbox{$\cdot$}}
\put(33.7,10.5){\mbox{$\cdot$}} \put(33.6,11){\mbox{$\cdot$}}
\put(33.5,11.5){\mbox{$\cdot$}} \put(33.4,12){\mbox{$\cdot$}}
\put(33.3,12.5){\mbox{$\cdot$}} \put(33.2,13){\mbox{$\cdot$}}
\put(33.1,13.5){\mbox{$\cdot$}} \put(33,14){\mbox{$\cdot$}}
\put(32.9,14.5){\mbox{$\cdot$}}

\put(35.9,0.5){\mbox{$\cdot$}} \put(36,1){\mbox{$\cdot$}}
\put(36.1,1.5){\mbox{$\cdot$}} \put(36.2,2){\mbox{$\cdot$}}
\put(36.3,2.5){\mbox{$\cdot$}} \put(36.4,3){\mbox{$\cdot$}}
\put(36.5,3.5){\mbox{$\cdot$}} \put(36.6,4){\mbox{$\cdot$}}
\put(36.7,4.5){\mbox{$\cdot$}} \put(36.8,5){\mbox{$\cdot$}}
\put(36.9,5.5){\mbox{$\cdot$}} \put(37,6){\mbox{$\cdot$}}
\put(37.1,6.5){\mbox{$\cdot$}} \put(37.2,7){\mbox{$\cdot$}}
\put(37.3,7.5){\mbox{$\cdot$}} \put(37.4,8){\mbox{$\cdot$}}
\put(37.5,8.5){\mbox{$\cdot$}} \put(37.6,9){\mbox{$\cdot$}}
\put(37.7,9.5){\mbox{$\cdot$}} \put(37.8,10){\mbox{$\cdot$}}
\put(37.9,10.5){\mbox{$\cdot$}} \put(38,11){\mbox{$\cdot$}}
\put(38.1,11.5){\mbox{$\cdot$}} \put(38.2,12){\mbox{$\cdot$}}
\put(38.3,12.5){\mbox{$\cdot$}} \put(38.4,13){\mbox{$\cdot$}}
\put(38.5,13.5){\mbox{$\cdot$}} \put(38.6,14){\mbox{$\cdot$}}
\put(38.7,14.5){\mbox{$\cdot$}}

\put(31,1){\mbox{$\Sigma_{\,4}$}}


\put(44.77,15){$_{\bullet}$} \put(56.77,15){$_{\bullet}$}
\put(50.77,14){$_{\bullet}$} \put(49.77,12){$_{\bullet}$}
\put(51.77,12){$_{\bullet}$} \put(50.77,9){$_{\bullet}$}
\put(50.77,0){$_{\bullet}$}

\put(49.2,11){\small 1} \put(52.2,11){\small 2} \put(52,13){\small
3} \put(56,15.5){$4$} \put(45,15.5){$5$} \put(51.8,8.2){$6$}
  \put(52,0){$7$}

\thicklines

\put(51,9){\line(-1,1){6}} \put(51,9){\line(1,1){6}}
\put(45,15){\line(1,0){12}}

\thinlines

\put(51,9){\line(-1,3){1}} \put(51,9){\line(1,3){1}}
\put(51,0){\line(-2,5){6}} \put(51,0){\line(2,5){6}}
\put(51,0){\line(0,1){9}} \put(50,12){\line(1,2){1}}
\put(50,12){\line(1,0){2}} \put(52,12){\line(-1,2){1}}
\put(45,15){\line(6,-1){6}} \put(45,15){\line(5,-3){5}}
\put(57,15){\line(-6,-1){6}} \put(57,15){\line(-5,-3){5}}

\put(46,1){\mbox{$\Sigma_{\,5}$}}

\end{picture}

\noindent $\Sigma_{\,1}, \dots, \Sigma_{\,5}$ {\rm are
combinatorial spheres. $\RR P^{\,2}_6$ triangulates the real
projective plane. $\Upsilon_1$, $\Upsilon_2$ are the smallest
examples of weak pseudomanifolds which are not pseudomanifolds. }
\end{eg}

The following results (which we need later) follow from the
classification of all 2-dimensional weak pseudomanifolds on $\leq
7$ vertices (e.g., see \cite{bd2, d}).

\begin{prop}$\!\!${\bf .} \label{6vertex}
Let $K$ be an $n$-vertex $2$-dimensional weak pseudomanifold. If
$n\leq 6$ then $K$ is isomorphic to $S^{\,2}_4$, $S^{\,1}_3\ast
S^{\,0}_2$, $S^{\,0}_2\ast S^{\,0}_2\ast S^{\,0}_2$, $\RR
P^{\,2}_6$ or $\Sigma_{\,1}$ above.
\end{prop}

\begin{prop}$\!\!${\bf .} \label{7vertex}
Let $K$ be a $7$-vertex $2$-dimensional weak pseudomanifold. If
the number of facets of $K$ is $\leq 10$ then $K$ is isomorphic to
$S^{\,1}_5 \ast S^{\,0}_2$, $\Sigma_{\,2}, \dots, \Sigma_{\,5}$,
$\Upsilon_1$ or $\Upsilon_2$ above.
\end{prop}

Let $X$ be a pure simplicial complex of dimension $d\geq 1$. Let
$A$ be a set of size $d+2$ such that $A$ contains at least one and
at most $d+1$ facets of $X$. (It follows that all except at most
one element of $A$ are vertices of $X$.) Define the pure
$d$-dimensional simplicial complex $\kappa_A(X)$ as follows. The
facets of $\kappa_A(X)$ are (i) the facets of $X$ not contained in
$A$ and (ii) the $(d+1)$-subsets of $A$ which are not facets of
$X$. $\kappa_A$ is said to be a {\em generalized bistellar move}.
Clearly $\kappa_A(\kappa_A(X)) = X$.  Let $\beta = \{x\in A \, :
\, A\setminus\{x\}\in X\}$ and $\alpha = A\setminus\beta$. Then
$\alpha\in X$ and $\beta\in \kappa_A(X)$. The set $\beta$ is
called the {\em core} of $A$. If $\alpha$ is an $i$-simplex of $X$
then $\kappa_A$ is also called a generalized bistellar {\em
$i$-move}. Observe that if $d$ is even and $\kappa_A$ is a
generalized bistellar $(d/2)$-move then $f_d(\kappa_A(X)) =
f_d(X)$.

Now suppose $X$ is a weak pseudomanifold, and $A$, $\alpha$ and
$\beta$ are as above. Notice that (a) either $\alpha$ is a
$d$-simplex in $X$ or $V({\rm Lk}_X(\alpha )) \supseteq \beta$ and
(b) if $\beta\in X$ then ${\rm Lk}_{\kappa_A(X)}(\beta) = {\rm
Lk}_{X}(\beta)\cup S^{\,i-1}_{i+1}(\alpha)\neq S^{\,i-1}_{i+
1}(\alpha)$ (and therefore $\kappa_A(X)$ is not a combinatorial
manifold even if $X$ is so). We shall say that $\kappa_A$ is a
{\em bistellar} move if ({\sf bs1}) $\beta \not \in X$ and ({\sf
bs2}) either $\alpha$ is a $d$-simplex in $X$ or $V({\rm
Lk}_X(\alpha)) =\beta$ (and hence ${\rm Lk}_X(\alpha)$ is the
standard sphere on the vertex set $\beta$). If $1\leq i\leq d-1$
then a bistellar $i$-move is called a {\em proper} bistellar move.
Observe that if $X$ is a combinatorial $d$-manifold then ({\sf
bs2}) holds for any $(d+2)$-subset $A$. If a generalized bistellar
move is not a bistellar move then it is called {\em singular}.

Two weak pseudomanifolds are called {\em bistellar equivalent} if
there exists a finite sequence of bistellar moves leading from one
to the other. Let $\kappa_A$ be a bistellar move on $X$. If $X_1$
is obtained from $X$ by starring (\cite{bd2}) a new vertex in
$\alpha$ and $X_2$ is obtained from $\kappa_A(X)$ by starring a
new vertex in $\beta$ then $X_1$ and $X_2$ are isomorphic. Thus if
$X$ and $Y$ are bistellar equivalent then $X\approx Y$. In
\cite{p}, Pachner proved the following\,: {\em Two combinatorial
manifolds are bistellar equivalent if and only if they are
combinatorially equivalent}.


\begin{eg}$\!\!${\bf .} \label{e2} {\rm Let the notations be as in
Example \ref{e1}.
\begin{enumerate}
   \item[$(a)$]  Let $A= \{1, 2, 5, 6\}\subset V(\RR P^{\,2}_6)$.
   Put $R = \kappa_A(\RR P^{\,2}_6)$. Then $R$ is not a weak
   pseudomanifold.  Observe that ({\sf bs1}) is not satisfied
   here and hence $\kappa_A$ is a singular bistellar move.
   Note that the automorphism group $A_5$ of $\RR P^{\,2}_6$ is
   transitive on the $4$-subsets of its vertex set. In consequence,
   all singular bistellar 1-moves on $\RR P^{\,2}_6$ yield
   isomorphic simplicial complexes.
     \item[$(b)$] Let $B=\{2, 3, 6, 7\}\subseteq V(\Sigma_{\,2})$. Then
     $\kappa_B(\Sigma_{\,2})$ is the union of two spheres with one
     common edge $67$. Here ({\sf bs1}) is not satisfied.
     \item[$(c)$] Let $C=\{1, 2, 3, 6\}\subseteq V(\Upsilon_1)$. Then
     $\kappa_C(\Upsilon_1)=\Upsilon_2$. Here also ({\sf bs1}) is not
     satisfied and $\kappa_C(\Upsilon_1) \not\approx \Upsilon_1$ but
     $\kappa_C(\Upsilon_1)$ is a weak pseudomanifold.
     \item[$(d)$] Let $D = \{1, 2, 3, 6\}\subseteq V(\Upsilon_2)$.
     Then $\kappa_D(\Upsilon_2)=\Upsilon_1$. Here ({\sf bs2}) is not
     satisfied.
     \item[$(e)$] If $E=\{2, 3, 4, 6\}\subseteq
     V(\Sigma_{\,4})$ then
     $\kappa_E(\Sigma_4)$ is a $7$-vertex pseudomanifold with 12
     facets. In this case, ({\sf bs1}) is not satisfied.
     \item[$(f)$] Let $F= \{2, 3, 4, 6\}\subseteq V(\Sigma_{\,2})$. Then
     $\kappa_F$ is a bistellar move and $\kappa_F(\Sigma_{\,2})=
     \Sigma_{\,3}$.
     \end{enumerate}
     }
\end{eg}

Let $L\subseteq K$ be simplicial complexes. The {\em simplicial
neighbourhood} of $L$ in $K$ is the subcomplex $N(L, K)$ of $K$
whose maximal simplices are those maximal simplices of $K$ which
intersect $V(L)$. Clearly, $N(L, K)$ is the smallest subcomplex of
$K$ whose geometric carrier is a topological neighbourhood of
$|L|$ in $|K|$. The induced subcomplex $C(L, K)$ on the vertex-set
$V(K)\setminus V(L)$ is called the {\em simplicial complement} of
$L$ in $K$.

Suppose $P^{\,\prime}\subseteq P$ are polyhedra and $P =
P^{\,\prime} \cup B$, where $B$ is a   pl $k$-ball (for some
$k\geq 1$). If $P^{\,\prime} \cap B$ is a   pl $(k-1)$-ball then
we say that there is an {\em elementary collapse} of $P$ to
$P^{\,\prime}$. We say that $P$ collapses to $Q$ and write $P\coll
Q$ if there exists a sequence $P = P_0, P_1, \dots, P_n = Q$ of
polyhedra such that there is an elementary collapse of $P_{i-1}$
to $P_{i}$ for $1\leq i\leq n$. If $Q$ is a point we say that $P$
is collapsible and write $P\coll 0$. For two simplicial complexes
$K$ and $L$, if $K\scoll L$  then clearly $|K| \coll |L|$. A {\em
regular neighbourhood} of a polyhedron $P$ in a pl $d$-manifold
$M$ is a $d$-dimensional submanifold $W$ with boundary such that
$W\coll P$ and $W$ is a neighbourhood of $P$ in $M$. The following
is a direct consequence of the Simplicial Neighbourhood Theorem
(\cite[Theorem 3.11]{rs}).

\begin{prop}$\!\!${\bf .} \label{snt} Let $K$ be a combinatorial
$d$-manifold with boundary. Suppose $\partial K$ is an induced
subcomplex of $K$. Let $L$ be the simplicial complement of
$\partial K$ in $K$. Then $|K| \coll |L|$.
\end{prop}

\noindent {\bf Proof.} Let $M$ be a pl d-manifold such that $|K|$
is in the interior of $M$ (we can always find such $M$, e.g., one
such $M$ can be obtained from $|K|\sqcup (|\partial K|\times [0,
1])$ by identifying $(x, 0)$ with $x\in |\partial K|$).

Since $L = C(\partial K, K)$, $|L| \subseteq |K|\setminus
|\partial K|$ and hence $|K|$ is a neighbourhood of $|L|$ in ${\rm
int}(M)$. Again, since $L$ is the simplicial complement of
$\partial K$ in $K$ and $\partial K$ is an induced subcomplex of
$K$, $C(L, K) = \partial K$. Finally, since $\partial K$ is an
induced subcomplex of dimension $d-1$, each $d$-simplex of $K$
intersects $V(L)$. This implies that $N(L, K) = K$.

Let $P = |L|$, $A= |K|$ and $J = \partial K$. Then $\partial A =
|\partial K|$ and $\bdN(L, K) := N(L, K)\cap C(L, K) = J$. Thus
(i) $P$ is a compact polyhedron in the interior of the pl manifold
$M$, (ii) $A$ is a neighbourhood of $P$ in ${\rm int}(M)$, (iii)
$A$ is a compact pl manifold with boundary and (iv) $(K, L, J)$
are triangulations of $(A, P, \partial A)$ where $L$ is an induced
subcomplex of $K$, $K=N(L, K)$ and $J= \bdN(L, K)$. Then, by the
Simplicial Neighbourhood Theorem, $A$ is a regular neighbourhood
of $P$. Hence $A\coll P$. \hfill $\Box$

\bigskip

We need the following well-known results (see \cite[Lemma 1.10,
Corollaries 3.13, 3.28]{rs}) later.

\begin{prop}$\!\!${\bf .} \label{ks1.10} Let $B$, $D$ be pl
$d$-balls and $h\colon \partial B\to \partial D$ a pl
homeomorphism. Then $h$ extends to a pl homeomorphism $h_1\colon
B\to D$.
\end{prop}

\begin{prop}$\!\!${\bf .} \label{ks3.13} Let $S$ be a pl
$d$-sphere. If $B\subseteq S$ is a pl $d$-ball then the closure of
$S\setminus B$ is a pl $d$-ball.
\end{prop}

\begin{prop}$\!\!${\bf .} \label{ks3.28} A collapsible pl
manifold with boundary is a pl ball.
\end{prop}

\noindent {\bf Question\,.} Is it true that under the hypothesis
of Proposition \ref{snt}, we have $K\scoll L$\,?

\section{$\ZZ_2$-acyclic simplicial complexes.}

In this section we prove Theorem \ref{t1}.

\begin{lemma}$\!\!${\bf .} \label{l3.1}
Let $X$ be a $7$-vertex simplicial complex. Suppose $(a)$ $X$ is
$\ZZ_2$-acyclic, $(b)$ $X$ is not collapsible, and $(c)$ $X$ is
minimal subject to $(a)$ and $(b)$ $($i.e., $X$ has no proper
subcomplex satisfying $(a)$ and $(b))$. Then $X$ is pure of
dimension $d=2$ or $3$ and each $(d-1)$-face of $X$ occurs in at
least two facets.
\end{lemma}

\noindent {\bf Proof.} Notice that, because of the minimality
assumption, $X$ has no free face. Clearly, $\dim(X)\leq 5$, since
otherwise $X$ is a combinatorial ball. Suppose $\dim(X) = 5$. By
minimality, each 4-face of $X$ is in 0 or $\geq 2$ facets. Since
$X$ has 7 vertices, it follows that each 4-face is in 0 or 2
facets. Therefore the pure part $Y$ of $X$ is a 7-vertex
5-dimensional weak pseudomanifold and hence $Y=S^{\,5}_7\subseteq
X$. Then $H_5(X, \ZZ_2)\not= 0$, a contradiction. Thus
$\dim(X)\leq 4$.

Suppose, if possible, $\dim(X)=4$. Let $Y$ be the pure part of
$X$. Then, each 3-face of $Y$ occurs in at least two facets. If
$\#(V(Y))\leq 6$, then $Y=S^{\,4}_6$ and hence $H_4(X, \ZZ_2)\not
= 0$, a contradiction. Thus $V(Y) = V(X)$ has size 7. Define a
binary relation $\sim$ on $V(Y)$ by $y_1\sim y_2$ if $V(Y)
\setminus\{y_1, y_2\}$ is not a facet of $Y$. Since each 3-face of
$Y$ is in at least two facets, it follows that $\sim$ is an
equivalence relation with at least two equivalence classes.
Therefore either there is an equivalence class $W$ of size 6 or
else we can write $V(Y)= V_1\sqcup V_2$, where $V_1$, $V_2$ are
unions of $\sim$-classes and $\#(V_1)\geq 2$,  $\#(V_2)\geq 2$. In
consequence $Y$ (and hence $X$) contains a 4-sphere as a
subcomplex\,: the standard sphere on $W$ or the join of the
standard spheres on $V_1$ and $V_2$. Therefore $H_4(X, \ZZ_2)\not
= 0$, a contradiction. Thus $\dim(X)\leq 3$.

If $\dim(X)=1$ then $X$ is a $\ZZ_2$-acyclic connected graph and
hence is a tree. But any tree has end vertices and hence is
collapsible, a contradiction. So, $\dim(X) = 2$ or 3.

Since $\widetilde{H}_0(X, \ZZ_2) = 0$, $X$ is connected. Since $X$
has no free vertex, it follows that each vertex of $X$ is in at
least two edges.

Next we show that $X$ has no maximal edge. Suppose, on the
contrary, $X$ has a maximal edge $e$. Then $Y := X\setminus \{e\}$
is a subcomplex of $X$. We claim that $Y$ is disconnected. If not,
then there is a subcomplex $K = S^{\,1}_n$ of $X$ containing the
edge $e$. The formal sum of the edges in $K$ is an 1-cycle over
$\ZZ_2$ which is not a boundary since it involves the maximal edge
$e$. Hence $H_1(X, \ZZ_2) \not = 0$, a contradiction. So, $Y$ is
disconnected. Since each vertex of $X$ is in at least two edges,
it follows that each component of $Y$ has $\geq 3$ vertices. Since
$X$ has seven vertices, it follows that some component of $Y$ has
exactly three vertices and contains an $S^{\,1}_3$. If these three
vertices span a 2-face then its edges are free in $X$,
contradicting minimality. In the remaining case $X$ has an induced
$S^{\,1}_3$ whose edges are maximal, contradicting
$\ZZ_2$-acyclicity of $X$.

In case $\dim(X)=2$, this shows that $X$ is pure. In case
$\dim(X)=3$, we proceed to show that $X$ has no maximal 2-face,
proving that it is pure in that case too.

Suppose, on the contrary, that $\dim(X)=3$ and $X$ has a maximal
2-face $\Delta = abc$. Let's say that an edge of $X$ is {\em good}
if it is in a tetrahedron of $X$, and call it {\em bad} otherwise.
First suppose that all three edges in $\Delta$ are good. Since $X$
has no free triangle, each vertex in the link of an edge has
degree 0 or $\geq 2$ and hence there are at least three vertices
of degree $\geq 2$ in the link of a good edge. Since $\Delta$ is
maximal, it follows that the link of each of the three edges in
$\Delta$ has $\geq 3$ vertices outside $\Delta$. Since, there are
only four vertices outside $\Delta$, it follows from the
pigeonhole principle that there is a common vertex $x$ outside
$\Delta$ which occurs in the link of all three edges in $\Delta$.
Hence $S^{\,2}_4(\Delta\cup\{x\})$ is a subcomplex of $X$. The sum
of the four triangles in this $S^{\,2}_4$ is a 2-cycle (with
$\ZZ_2$ coefficients) which can not be the boundary of a 3-chain
since one of these triangles is maximal. Therefore $H_2(X,
\ZZ_2)\not = 0$, a contradiction. Thus $\Delta$ contains at least
one bad edge.

We claim that $\Delta$ can't have more than one bad edges.
Suppose, on the contrary, that $ab$ and $ac$ are bad edges in $X$.
Notice that (arguing as in the proof of the case $\dim(X)=4$), if
a 3-dimensional simplicial complex on $\leq 6$ vertices has $\geq
2$ tetrahedra through each triangle then it contains a
combinatorial $S^{\,3}$. Therefore the pure part $Y$ of $X$ must
have seven vertices. In particular $a\in Y$. Since $ab$ and $ac$
are bad edges, $b, c\not\in {\rm Lk}_Y(a)$ and hence
$\deg_Y(a)\leq 4$. Therefore ${\rm Lk}_Y(a)=S^{\,2}_4$. Hence we
can apply an improper bistellar move to $Y$ to remove the vertex
$a$, yielding a 6-vertex 3-dimensional simplicial complex
$\widetilde{Y}$ with $\geq 2$ tetrahedra through each triangle.
Hence $\widetilde{Y}$ has an $S^{\,3}$ as a subcomplex, so that
$H_3(Y, \ZZ_2)= H_3(\widetilde{Y}, \ZZ_2)\not = 0$. Therefore
$H_3(X, \ZZ_2) \not = 0$, a contradiction. Thus $\Delta$ contains
exactly one bad edge, say $ab$. Hence $ac$ and $bc$ are good
edges.

Since $X$ has no free edge, there is a second triangle, say $abd$,
through $ab$. Since $ab$ is a bad edge, $abd$ is maximal. By the
above argument, $ad$ and $bd$ are good edges. If both $acd$ and
$bcd$ are triangles of $X$ then $X$ has $S^{\,2}_4(a, b, c, d)$ as
a subcomplex, and at least one of the triangles of this
$S^{\,2}_4$ is maximal in $X$, yielding the contradiction $H_2(X,
\ZZ_2)\not = 0$ as before. Therefore, without loss of generality,
we may assume $bcd\not\in X$. Note that $a$ is an isolated vertex
in ${\rm Lk}_X(bc)$ and $d$ does not occur in ${\rm Lk}_X(bc)$.
Since $bc$ is a good edge, it follows that all three vertices
outside $\{a, b, c, d\}$ (say $x$, $y$ and $z$) occur in ${\rm
Lk}_X(bc)$. Similarly, $x, y, z\in {\rm Lk}_X(bd)$. Again, the
good edges $ac$ and $ad$ have at most one non-isolated vertex from
$\{a, b, c, d\}$ in their links, hence each of them has at least
two of $x$, $y$, $z$ in their links. Therefore, there is one
vertex, say $x$, which occurs in the link of all the four edges
$ac$, $bc$, $ad$, $bd$. Hence $S^{\,0}_2(c, d)\ast S^{\,1}_3(a, b,
x)$ is a subcomplex of $X$. Since one of the triangles in this
2-sphere is maximal, it follows that $H_2(X, \ZZ_2)\not = 0$, a
contradiction. Thus $X$ has no maximal triangles nor maximal
edges, so $X$ is pure.

Finally, the last assertion follows from purity and minimality of
$X$. \hfill $\Box$

\begin{lemma}$\!\!${\bf .} \label{l3.2}
Let $X$ be a $7$-vertex $2$-dimensional $\ZZ_2$-acyclic simplicial
complex. Then $X$ is collapsible.
\end{lemma}

\noindent {\bf Proof.} Let $X$ be a minimal counter example. Let
$f_i$, $0\leq i\leq 2$, be the number of $i$-faces in $X$. Since
$X$ is $\ZZ_2$-acyclic, $\chi(X) = 1$. Thus, $f_0=7$ and
$f_1=f_2+6$.

For $i\geq 0$, let $e_i$ be the number of edges of degree $i$ in
$X$. By Lemma \ref{l3.1}, $e_i = 0$ for $i\leq 1$. Two-way
counting yields
$$
\sum_{i=2}^5e_i=f_1=f_2+6, ~~ \sum_{i=2}^5i\hspace{0.4mm}e_i=3f_2.
$$
Hence
\begin{equation} \label{eq3.2a}
e_3+3e_5\leq e_3+2e_4+3e_5=f_2-12.
\end{equation}
Let's say that an edge of $X$ is {\em odd} (respectively {\em
even}) if it lies in an odd (respectively even) number of
triangles. Note that each graph has an even number of vertices of
odd degree. Applying this trivial observation to the vertex links
of $X$, we conclude that each vertex of $X$ is in an even number
of odd edges. Thus the total number  $e_3+e_5$ of odd edges is $=
0$ or $\geq 3$. If there is no odd edge then the sum of all the
triangles gives a non-zero element of $H_2(X, \ZZ_2)$, a
contradiction. So, $e_3+e_5\geq 3$. Combining this with
(\ref{eq3.2a}), we get $f_2\geq 15$ and hence $f_1\geq 21 = {7
\choose 2}$. Hence $f_1=21$, $f_2=15$, $e_3=3$, $e_4=e_5=0$.

Since each vertex is in an even number of odd edges, it follows
that the three odd edges form a triangle $\Delta$, which may or
may not be in $X$.

If $\Delta$ is in $X$, then the sum of the remaining triangles
gives a non-zero element of $H_2(X, \ZZ_2)$, a contradiction. If
$\Delta$ is not in $X$ then (as each of the three edges in
$\Delta$ has three vertices in its link and there are four
vertices outside $\Delta$) by the pigeonhole principle there is a
vertex $x\not\in \Delta$ such that $x$ occurs in the link of each
of the three edges in $\Delta$. Then the sum of all the triangles
excepting the three triangles in $\Delta\cup\{x\}$ gives a
non-zero element of $H_2(X, \ZZ_2)$, a contradiction. \hfill
$\Box$

\begin{lemma}$\!\!${\bf .} \label{l3.3}
Let $U$ be a $2$-dimensional pure simplicial complex on $\leq 7$
vertices. Suppose the number of triangles in $U$ is $\leq 10$ and
each edge of $U$ is in an even number of triangles. Then either
$U$ is the union of two combinatorial spheres $($on $4$ or $5$
vertices$)$ with no common triangle, or $U$ is isomorphic to one
of $S^{\,2}_4$, $S^{\,1}_3\ast S^{\,0}_2$, $S^{\,0}_2\ast
S^{\,0}_2 \ast S^{\,0}_2$, $S^{\,1}_5 \ast S^{\,0}_2$, $\RR
P^{\,2}_6$, $\Sigma_{\,1}, \dots, \Sigma_{\,5}$ or $R$ $($of
Example $\ref{e1}$ and Example $\ref{e2}$ $(a))$.
\end{lemma}

\noindent {\bf Proof.} Let ${\cal S}$ be the list of simplicial
complexes in the statement of this lemma. We find by inspection
that ${\cal S}$ is closed under generalized bistellar 1-moves.

If $f_0(U) \leq 5$ then $U$ is a weak pseudomanifold and hence, by
Proposition \ref{6vertex}, $U\in {\cal S}$. So assume $f_0(U) = 6$
or 7. The proof is by induction on the number $n(U)$ of degree 4
edges in $U$. If $n(U)=0$ then $U$ is a weak pseudomanifold and
hence, by Propositions \ref{6vertex} and \ref{7vertex}, $U\in
{\cal S}$. So let $n(U)>0$ and suppose that we have the result for
all smaller values of $n(U)$.

By the assumption, all the edges of $U$ are of degree 2 or 4.
Therefore, a two-way counting yields $4n(U)+2(f_1(U)-n(U)) =
3f_2(U)\leq 30$. Thus, $n(U)+f_1(U)\leq 15$. Therefore,
\begin{equation} \label{eq3.3a}
f_1(U)< 15,
\end{equation}
showing that $U$ has at least one non-edge. Fix an edge $ab$ of
degree 4 in $U$. Let $W$ be the link of $ab$. If each pair of
vertices in $W$ formed an edge in $U$ then $f_1(U)$ would be $\geq
15$, contradicting (\ref{eq3.3a}). So, there exist $c, d\in W$
such that $cd$ is a non-edge in $U$.

Let $A = \{a, b, c, d\}$. Then $\kappa_A$ is a generalized
bistellar 1-move and hence $\kappa_A(U)$ also satisfies the
hypothesis of the lemma, and $n(\kappa_A(U)) = n(U) -1$.
Therefore, by the induction hypothesis, $\kappa_A(U) \in {\cal
S}$. Since ${\cal S}$ is closed under generalized bistellar
1-moves, $U = \kappa_A(\kappa_A(U))\in {\cal S}$. \hfill $\Box$

\begin{lemma}$\!\!${\bf .} \label{l3.4}
Let $X$ be a $7$-vertex $3$-dimensional simplicial complex.
Suppose $(a)$ $X$ is $\ZZ_2$-acyclic, $(b)$ $X$ is not
collapsible, and $(c)$ $X$ is minimal subject to $(a)$ and $(b)$.
Then the $f$-vector of $X$ is $(7, 20, 30, 16)$, $(7, 21, 32,
17)$, $(7, 21, 33, 18)$, $(7, 21, 34, 19)$ or  $(7, 21, 35, 20)$.
\end{lemma}

\noindent {\bf Proof.} For $0\leq i\leq 3$, let $f_i$ be the
number of $i$-faces of $X$. For $i\geq 0$, let $t_i$ be the number
of triangles of degree $i$ in $X$. By Lemma \ref{l3.1}, we have
$t_i=0$ for $i\leq 1$. Two way counting yields
$$
\sum_{i=2}^4 t_i = f_2, ~~ \sum_{i=2}^4 i\hspace{0.4mm}t_i = 4 f_3
$$
and hence
\begin{equation} \label{eq3.4a}
t_3\leq t_3 + 2 t_4 = 4 f_3 - 2f_2.
\end{equation}

Say that a triangle of $X$ is {\em odd} (respectively {\em even})
if it is in an odd (respectively  even) number of tetrahedra of
$X$. By the same argument as in Lemma \ref{l3.2}, each edge is in
an even number of odd triangles, so that the number $t_3$ of odd
triangles is 0 or $\geq 4$.

If there is no odd triangle then the sum of all the tetrahedra
gives a non-zero element of $H_3(X, \ZZ_2)$, a contradiction. So,
$t_3\geq 4$. Combining this with (\ref{eq3.4a}) we get
\begin{equation} \label{eq3.4b}
 2 f_3 - f_2 \geq 2.
\end{equation}

Since $X$ is $\ZZ_2$-acyclic, by a result of Stanley (\cite{st}),
$X$ has a 2-dimensional subcomplex $Y$ such that the $f$-vector of
$X$ equals the $f$-vector of a cone over $Y$. (In \cite{st}, the
author uses the vanishing of the reduced cohomology groups as his
definition of acyclicity, while we have used the homology
definition. However, since the coefficient ring  used is a field,
these two definitions coincide.) Let $(g_0, g_1, g_2)$ be the
$f$-vector of $Y$. Thus, $g_0=6$ and
\begin{equation} \label{eq3.4c}
f_1 = g_1 + 6, ~~ f_2 = g_1 + g_2, ~~ f_3 = g_2.
\end{equation}
Hence (\ref{eq3.4b}) yields
\begin{equation} \label{eq3.4d}
g_2\geq g_1 +2.
\end{equation}
Let $m={6\choose 2}-g_1$, $n={6\choose 3}-g_2$ be the number of
non-edges and non-triangles of $Y$, respectively. Since each
non-edge is in exactly four non-triangles and any two non-edges
are shared by at most one non-triangle, we have $n\geq
4m-{m\choose 2}$. Also, from (\ref{eq3.4d}) we get $n\leq m+3$.
Hence $m+3\geq 4m-{m\choose 2}$ or $(m-1)(m-6)\geq 0$. So, either
$m\leq 1$ or $m\geq 6$.

First suppose $m\geq 6$, i.e., $g_1\leq 9$. If each edge of $Y$
was in $\leq 3$ triangles then we would have $g_2\leq g_1$,
contradicting (\ref{eq3.4d}). So, there is an edge of $Y$
contained in four triangles, together covering all the nine edges
of $Y$. But, apart from the four triangles already seen, no three
of these nine edges form a triangle of $Y$. Thus $g_2=4$, $g_1=9$
-- contradicting (\ref{eq3.4d}). So, $m\leq 1$, i.e., $g_1=14$ or
15.

If $g_1=14$ then the four triangles through the missing edge are
missing from $Y$, so that $g_2\leq 16$. Thus, by (\ref{eq3.4d}),
$(g_1, g_2)= (14, 16)$, $(15, 17)$, $(15, 18)$, $(15, 19)$ or
$(15, 20)$. The lemma now follows from (\ref{eq3.4c}). \hfill
$\Box$

\begin{lemma}$\!\!${\bf .} \label{l3.5}
Let $X$ be a $7$-vertex $3$-dimensional $\ZZ_2$-acyclic simplicial
complex. Then $X$ is collapsible.
\end{lemma}

\noindent {\bf Proof.} Let $X$ be a minimal counter example. As
before, each edge is in an even number of odd triangles. Let
$f_i$'s and $t_j$'s be as in the proof of Lemma \ref{l3.4}. Then,
by Lemma \ref{l3.4}, $t_3+2t_4=4f_3-2f_2\leq 10$ and hence the
number $t_3$ of odd triangles is $\leq 10$.

Let $U$ denote the pure 2-dimensional simplicial complex whose
facets are the odd triangles of $X$.  Then each edge of $U$ is in
an even number of triangles of $U$. Therefore, by Lemma
\ref{l3.3}, we get the following cases\,:

\medskip

\noindent {\sf Case 1\,:} $U$ is the union of two combinatorial
spheres with no common triangle (on 4 or 5 vertices), say on
vertex sets $A$ and $B$.

First suppose $\#(A)=\#(B)=4$. If both $A$ and $B$ are 3-faces in
$X$ then the pure simplicial complex $\widetilde{X}$ whose facets
are those of $X$ other than $A$, $B$ is a 3-dimensional weak
pseudomanifold. This implies that the sum of all the tetrahedra,
excepting $A$ and $B$, gives a non-zero element of $H_3(X,
\ZZ_2)$, a contradiction. So, without loss of generality $A \not
\in X$.

Since each of the four triangles inside $A$ is of degree 3 in $X$,
the three vertices (say $x$, $y$, $z$) outside $A$ occur in the
link of all the four triangles. Then the 3-sphere
$S^{\,2}_4(A)\ast S^{\,0}_2(x, y)$ occurs as a subcomplex of $X$,
forcing $H_3(X, \ZZ_2)\not=0$, a contradiction.

In the remaining case $\#(A)= 4$, $\#(B)=5$ (since $U$ has at most
10 triangles, the case $\#(A)=\#(B)=5$ does not arise). Write
$B=\{b_1, b_2, b_3, x, y\}$ and $U= S^{\,2}_4(A) \cup
(S^{\,1}_3(b_1, b_2, b_3)\ast S^{\,0}_2(x, y))$. As above, we must
have $A\in X$.

If both $b_1b_2b_3x$ and $b_1b_2b_3y$ are in $X$, then the sum of
the 3-faces other than $A$, $b_1b_2b_3x$ and $b_1b_2b_3y$ gives a
non-zero element of $H_3(X, \ZZ_2)$, a contradiction. So, without
loss of generality, $b_1b_2b_3x \not\in X$. Since the triangles of
$S^{\,1}_3(b_1, b_2, b_3)\ast S^{\,0}_2(x, y))$ are degree 3
triangles in $X$, it follows that $b_1b_2xy$, $b_1b_3xy$,
$b_2b_3xy\in X$. Then the sum of the tetrahedra other than $A$ and
these three tetrahedra gives a non-zero element of $H_3(X,
\ZZ_2)$, a contradiction.

\medskip

\noindent {\sf Case 2\,:}  $U=S^{\,2}_4$. We get a contradiction
as in Case 1.

\medskip

\noindent {\sf Case 3\,:} $U=S^{\,1}_3\ast S^{\,0}_2$. We get a
contradiction as in Case 1.

\medskip

\noindent {\bf Observation 1\,:} As $t_3\geq 8$ in the remaining
cases, we have $2f_3 - f_2 \geq 4$ and hence only the following
two possibilities survive for the $f$-vector of $X$\,: $(7, 21,
34, 19)$ and $(7, 21, 35, 20)$. Therefore $X$ has at most one
missing triangle and at most one triangle of degree 4, and these
two cases are exclusive. It follows that, if $x$ is a vertex not
covered by the odd triangles, then ${\rm Lk}_X(x)$ is a 6-vertex
2-dimensional neighbourly weak pseudomanifold. But, from
Proposition \ref{6vertex}, we see that $\RR P^{\,2}_6$ is the only
possibility. Thus, ${\rm Lk}_X(x)=\RR P^{\,2}_6$. This implies
that if $V_1\subseteq V(U)$ is a 3-set then exactly one of $V_1$
and $V(U)\setminus V_1$ is a simplex in ${\rm Lk}_X(x)$. In
particular, any two triangles in ${\rm Lk}_X(x)$ intersect.

\medskip

\noindent {\sf Case 4\,:} $U=S^{\,0}_2(a_1, a_2)\ast
S^{\,0}_2(b_1, b_2)\ast S^{\,0}_2(c_1, c_2)$. Then the odd
triangles of $X$ are $a_ib_jc_k$, $1\leq i, j, k\leq 2$. If
$\{a_1a_2b_jc_k \, : \, 1\leq j, k\leq 2\}\subseteq X$, then the
sum of the remaining tetrahedra gives a non-zero element of
$H_3(X, \ZZ_2)$, a contradiction. So, without loss of generality,
$a_1a_2b_1c_1\not\in X$. As $a_1b_1c_1$, $a_2b_1c_1$ are degree 3
triangles, it follows that $a_1b_1b_2c_1$, $a_2b_1b_2c_1\in X$. If
both $a_1b_1b_2c_2$ and $a_2b_1b_2c_2$ are in $X$ then
$X\supseteq\{a_ib_1b_2c_k \, : 1\leq i, k\leq 2\}$, hence we get a
contradiction as before. So, without loss of generality,
$a_2b_1b_2c_2\not\in X$.

Since $a_1a_2b_1c_1$, $a_2b_1b_2c_2\not\in X$ and $a_1b_1c_1$,
$a_2b_2c_2$ are degree 3 triangles, it follows that these two
disjoint triangles occur in the link of $x$. But this contradicts
Observation 1.

\medskip

\noindent {\sf Case 5\,:} $U=\Sigma_1$ of Example \ref{e1}. Thus,
the odd triangles are $125$, $126$, $156$, $235$, $236$, $345$,
$346$ and $456$. If $1256$, $3456\not\in X$ then, since $125$ and
$346$ are degree 3 triangles, they are disjoint triangles in ${\rm
Lk}_X(x)$, contradicting Observation 1. So, without loss of
generality, $1256\in X$.

If $3456\not\in X$ then, since $345$, $346$, $456$ are degree 3
triangles, $2345$, $2346$, $2456\in X$. Then the sum of all the
tetrahedra, excepting $1256$, $2345$, $2346$, $2456$, gives a
non-zero element of $H_3(X, \ZZ_2)$. So, $3456\in X$.

If $2356\in X$, then the sum of all the tetrahedra, excepting
$1256$, $2356$, $3456$, gives a non-zero element of $H_3(X,
\ZZ_2)$. Therefore $2356\not\in X$.

Since $235$ and $236$ are degree 3 triangles, $2345$, $2346\in X$.
First suppose that at least one of $1356$, $2456$ is in $X$.
Without loss, say $2456\in X$. Then the sum of all the tetrahedra,
excepting $1256$, $2456$, $2345$, $2346$, gives a non-zero element
of $H_3(X, \ZZ_2)$. Thus $1356$, $2456\not\in X$. Then, since
$156$, $456$ are degree 3 triangles, $156x$, $456x\in X$.

Since $2356, 2456\not\in X$, $x\in {\rm Lk}_X(256)$, i.e.,
$256x\in X$. Similarly, looking at $356$, we conclude that
$356x\in X$. Thus, $56x$ is a degree 4 triangle in $X$. But this
is not possible since,  by Observation 1, ${\rm Lk}_X(x)$ is $\RR
P^{\,2}_6$.

\medskip

\noindent {\bf Observation 2\,:} In the remaining cases, $t_3 =
10$ and hence the $f$-vector of $X$ is $(7, 21, 35, 20)$. In
consequence, $t_4=0$. Thus all triangles are of degree 2 or 3.
Since $f_3= {7\choose 3}$, each edge in $X$ has degree 5. Thus if
$e$ is an edge outside $U$ then the link of $e$ is a pentagon
($S^{\,1}_5$).

\medskip

\noindent {\sf Case 6\,:} $U=\RR P^{\,2}_6$. In this case, all the
4-sets of vertices not containing $x$ contain exactly two odd
triangles each. In particular, all the tetrahedra of $X$ not
containing $x$ contain exactly two odd triangles each. Trivially,
each tetrahedron through $x$ contains at most one odd triangle.
Thus, letting $\alpha_i$, $i\geq 0$, denote the number of
tetrahedra of $X$ containing exactly $i$ odd triangles, we have
$\alpha_2=20-10=10$ and $\alpha_0+\alpha_1=10$. But two way
counting yields $\alpha_1+2\alpha_2=10\times 3=30$. Hence
$\alpha_1=10$, $\alpha_0=0$. Thus $x$ occurs in the link of each
odd triangle and hence ${\rm Lk}_X(x)=U$. Therefore the 10
tetrahedra of $X$ not passing through $x$ add up to a non-zero
element of $H_3(X, \ZZ_2)$, a contradiction.

\medskip

\noindent {\sf Case 7\,:} $U = R$ of Example \ref{e2} $(a)$. Thus,
the odd triangles are $123$, $124, 125, 126, 135, 146$, $236$,
$245$, $345$ and $346$. We claim that ${\rm Lk}_X(12)\supseteq$
\setlength{\unitlength}{2mm}
\begin{picture}(5,2.7)(0,-0.1)

\put(0.7,0){$_{\bullet}$} \put(0.7,2){$_{\bullet}$}
\put(2.7,0){$_{\bullet}$} \put(2.7,2){$_{\bullet}$}

\thinlines

\put(1,0){\line(0,1){2}} \put(3,0){\line(0,1){2}}
\put(1,0){\line(1,0){2}} \put(1,2){\line(1,0){2}}

\put(3.7,1.5){\scriptsize $3$} \put(-0.2,-0.2){\scriptsize $4$}
\put(3.7,-0.2){\scriptsize $6$} \put(-0.2,1.5){\scriptsize $5$}
\end{picture}.
If, for instance, $1236\not\in X$ then,  as $123$, $126$, $236$
are degree 3 triangles, $x$ belongs to the link of each of these
triangles. Then ${\rm Lk}_X(2x)\supseteq$
\setlength{\unitlength}{2mm}
\begin{picture}(5,2.7)(0,-0.1)

\put(0.7,1){$_{\bullet}$} \put(2.7,0){$_{\bullet}$}
\put(2.7,2){$_{\bullet}$}

\thinlines

\put(3,0){\line(0,1){2}} \put(1,1){\line(2,1){2}}
\put(1,1){\line(2,-1){2}}

\put(3.7,1.6){\scriptsize $6$} \put(0.4,-0.3){\scriptsize $1$}
\put(3.7,-0.2){\scriptsize $3$}
\end{picture},
contradicting Observation 2. This proves the claim.

Since $3, 4, 5, 6$ are of degree 3 and $x$ is of degree 2 in ${\rm
Lk}_X(12)$, it follows that ${\rm Lk}_X(12)=$
\setlength{\unitlength}{2.5mm}
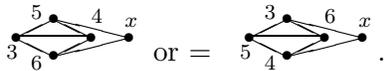
\begin{picture}(7.8,3.2)(0,-0.2)

\put(0.7,1){$_{\bullet}$} \put(2.7,0){$_{\bullet}$}
\put(2.7,2){$_{\bullet}$} \put(4.7,1){$_{\bullet}$}
\put(6.7,1){$_{\bullet}$}

\thinlines

\put(1,1){\line(2,1){2}} \put(1,1){\line(2,-1){2}}
\put(5,1){\line(-2,1){2}} \put(5,1){\line(-2,-1){2}}
\put(7,1){\line(-4,1){4}} \put(7,1){\line(-4,-1){4}}
\put(1,1.1){\line(1,0){4}}

\put(0.5,-0.1){\scriptsize {$3$}} \put(1.8,2){\scriptsize {$5$}}
\put(6.8,1.6){\scriptsize {$x$}} \put(1.8,-0.7){\scriptsize {$6$}}
\put(5,1.8){\scriptsize {$4$}}
\end{picture} or $=$
\setlength{\unitlength}{2.5mm}
\begin{picture}(7.8,3.2)(0,-0.2)

\put(0.7,1){$_{\bullet}$} \put(2.7,0){$_{\bullet}$}
\put(2.7,2){$_{\bullet}$} \put(4.7,1){$_{\bullet}$}
\put(6.7,1){$_{\bullet}$}

\thinlines

\put(1,1){\line(2,1){2}} \put(1,1){\line(2,-1){2}}
\put(5,1){\line(-2,1){2}} \put(5,1){\line(-2,-1){2}}
\put(7,1){\line(-4,1){4}} \put(7,1){\line(-4,-1){4}}
\put(1,1.1){\line(1,0){4}}

\put(0.5,-0.1){\scriptsize {$5$}} \put(1.8,2){\scriptsize {$3$}}
\put(6.8,1.6){\scriptsize {$x$}} \put(1.8,-0.7){\scriptsize {$4$}}
\put(5,1.8){\scriptsize {$6$}}
\end{picture}.

In the first case, $125, 126\in {\rm Lk}_X(x)$. Hence, by
Observation 1, $345$, $346\not\in {\rm Lk}_X(x)$. Since these two
are degree 3 triangles, it follows that ${\rm Lk}_X(345)= \{1, 2,
6\}$ and ${\rm Lk}_X(346)= \{1, 2, 5\}$. Since $1$, $2$ are of
degree 2 in ${\rm Lk}_X(34)$, this forces ${\rm Lk}_X(34)=$
\setlength{\unitlength}{2.5mm}
\begin{picture}(5.8,3.2)(0,-0.2)

\put(0.7,1){$_{\bullet}$} \put(2.7,0){$_{\bullet}$}
\put(2.7,2){$_{\bullet}$} \put(4.7,1){$_{\bullet}$}

\thinlines

\put(1,1){\line(2,1){2}} \put(1,1){\line(2,-1){2}}
\put(5,1){\line(-2,1){2}} \put(5,1){\line(-2,-1){2}}
\put(1,1.1){\line(1,0){4}}

\put(0.5,-0.1){\scriptsize {$5$}} \put(1.8,2){\scriptsize {$1$}}
\put(3.8,-0.6){\scriptsize {$2$}} \put(5,1.5){\scriptsize {$6$}}
\end{picture}
and hence $x\not\in {\rm Lk}_X(34)$. This is a contradiction since
$X$ is 3-neighbourly.

In the second case, $125, 126\not\in {\rm Lk}_X(x)$ and hence, by
Observation 1, $345, 346\in {\rm Lk}_X(x)$. That is, $5x$, $6x\in
{\rm Lk}_X(34)$. Also, as $34\not\in {\rm Lk}_X(12)$, we have
$12\not\in {\rm Lk}_X(34)$. Since $5, 6$ are of degree 3 and $1$,
$2$, $x$ are of degree 2 in ${\rm Lk}_X(34)$, it follows that
${\rm Lk}_X(34)=$ \setlength{\unitlength}{2.5mm}
\begin{picture}(7.8,3.2)(0,-0.2)

\put(0.7,1){$_{\bullet}$} \put(2.7,0){$_{\bullet}$}
\put(2.7,2){$_{\bullet}$} \put(4.7,1){$_{\bullet}$}
\put(6.7,1){$_{\bullet}$}

\thinlines

\put(1,1){\line(2,1){2}} \put(1,1){\line(2,-1){2}}
\put(5,1){\line(-2,1){2}} \put(5,1){\line(-2,-1){2}}
\put(7,1){\line(-4,1){4}} \put(7,1){\line(-4,-1){4}}

\put(0.5,-0.1){\scriptsize {$1$}} \put(1.8,2){\scriptsize {$5$}}
\put(6.8,1.5){\scriptsize {$x$}} \put(1.8,-0.7){\scriptsize {$6$}}
\put(3.3,0.6){\scriptsize {$2$}}
\end{picture}.
Hence $1345, 2345, 345x \in X$. Also, as $123$ is a degree 3
triangle and $1234\not\in X$, we have $1235\in X$. Thus
\setlength{\unitlength}{2mm}
\begin{picture}(6,3.1)(0,-0.1)

\put(0.7,1){$_{\bullet}$} \put(2.7,0){$_{\bullet}$}
\put(2.7,2){$_{\bullet}$} \put(4.7,2){$_{\bullet}$}

\thinlines

\put(3,2){\line(1,0){2}} \put(1,1){\line(2,1){2}}
\put(1,1){\line(2,-1){2}} \put(3,0){\line(0,1){2}}

\put(0.3,-0.3){\scriptsize {$2$}} \put(1.8,2){\scriptsize {$4$}}
\put(5.7,2){\scriptsize {$x$}} \put(3.4,0.4){\scriptsize {$1$}}
\end{picture}
$\subseteq {\rm Lk}_X(35)$. Since $1$, $4$ are of degree 3 while
$2$, $6$, $x$ are of degree 2 in this link, it follows that ${\rm
Lk}_X(35)=$ \setlength{\unitlength}{2mm}
\begin{picture}(6.6,3.1)(0,-0.1)

\put(0.7,1){$_{\bullet}$} \put(2.7,0){$_{\bullet}$}
\put(2.7,2){$_{\bullet}$} \put(4.7,0){$_{\bullet}$}
\put(4.7,2){$_{\bullet}$}

\thinlines

\put(3,0){\line(1,0){2}} \put(5,0){\line(0,1){2}}
\put(3,2){\line(1,0){2}} \put(1,1){\line(2,1){2}}
\put(1,1){\line(2,-1){2}} \put(3,0){\line(0,1){2}}

\put(0.3,-0.3){\scriptsize {$2$}} \put(1.8,2){\scriptsize {$4$}}
\put(5.7,2){\scriptsize {$x$}} \put(3.4,0.4){\scriptsize {$1$}}
\put(5.7,-0.1){\scriptsize {$6$}}
\end{picture}.
Hence $356x\in X$. Then
    \setlength{\unitlength}{2mm}
\begin{picture}(5,2.7)(0,-0.1)

\put(0.7,1){$_{\bullet}$} \put(2.7,0){$_{\bullet}$}
\put(2.7,2){$_{\bullet}$}

\thinlines

\put(3,0){\line(0,1){2}} \put(1,1){\line(2,1){2}}
\put(1,1){\line(2,-1){2}}

\put(3.7,1.6){\scriptsize $6$} \put(0.4,-0.3){\scriptsize $4$}
\put(3.7,-0.2){\scriptsize $5$}
\end{picture}
$\subseteq {\rm Lk}_X(3x)$, contradicting Observation 2.

\medskip

\noindent {\sf Claim\,:} In the remaining cases, if $F$ is a set
of four vertices of $U$ containing at least two odd triangles,
then either $F \in X$ or $F\subseteq V({\rm Lk}_U(x))$ for some
vertex $x$.

\smallskip

In these cases, $V(U)= V(X)$. If $F\not\in X$ contains two odd
triangles, then on the average, a vertex outside $F$ occurs in the
links (in $X$) of $\geq \frac{3\times 2 + 2\times 2}{3} > 3$ of
the four triangles inside $F$. Thus there is a vertex $x$ in the
link of all these triangles. If $F\not\subseteq V({\rm Lk}_U(x))$
for this $x$, then choose a vertex $y\in F$ such that $xy\not\in
U$. Then ${\rm Lk}_X(xy)\supseteq S^{\,1}_3(F\setminus\{y\})$,
contradicting Observation 2. This proves the claim.

\medskip

\noindent {\sf Case 8\,:} $U=S^{\,1}_5(\ZZ_5) \ast S^{\,0}_2(u,
v)$. In this case, the above claim implies that $X$ contains the
five tetrahedra $\{u, v, i, i+1\}$, $i\in \ZZ_5$. Then the sum of
the remaining fifteen tetrahedra gives a non-zero element of
$H_3(X, \ZZ_2)$, a contradiction.

\medskip

\noindent {\sf Case 9\,:} $U=\Sigma_{\,2}$ of Example \ref{e1}.
Thus, the odd triangles are $126$, $127$, $167$, $236$, $237$,
$346$, $347$, $456$, $457$ and $567$. By the above claim, $1267$,
$2367$, $3467$, $4567 \in X$. Then the sum of the remaining
sixteen tetrahedra gives a non-zero element of $H_3(X, \ZZ_2)$, a
contradiction.

\medskip

\noindent {\sf Case 10\,:} $U=\Sigma_{\,3}$ of Example \ref{e1}.
Thus, the odd triangles are $126$, $127$, $167$, $234$, $237$,
$246$, $347$, $456$, $457$ and $567$. By the claim, $1267$,
$2347$, $4567 \in X$.

If $2467\in X$ then the sum of all the tetrahedra, excepting
$1267$, $2347$, $4567$, $2467$, gives a non-zero element of
$H_3(X, \ZZ_2)$, a contradiction. So, $2467\not\in X$. Then, ${\rm
Lk}_X(246) = \{1, 3, 5\}$.

Since $\deg(247)=2$ and $2347\in X$, assume without loss of
generality, that $2457\in X$ and $1247\not\in X$. Then ${\rm
Lk}_X(127) = \{3, 5, 6\}$.

So, $2456$, $2457\in X$ and $\deg(245)=2$. Hence $2345\not\in X$.
Then ${\rm Lk}_X(234) = \{1, 6, 7\}$.

Now, $1234$, $1237\in X$ and $\deg(123) =2$. Therefore, $1236
\not\in X$. Then ${\rm Lk}_X(126) = \{4, 5, 7\}$. This implies
that \setlength{\unitlength}{2mm}
\begin{picture}(5,2.7)(0,-0.1)

\put(0.7,0){$_{\bullet}$} \put(0.7,2){$_{\bullet}$}
\put(2.7,0){$_{\bullet}$} \put(2.7,2){$_{\bullet}$}

\thinlines

\put(1,0){\line(0,1){2}} \put(3,0){\line(0,1){2}}
\put(1,0){\line(1,0){2}} \put(1,2){\line(1,0){2}}

\put(3.7,1.5){\scriptsize $1$} \put(-0.2,-0.2){\scriptsize $4$}
\put(3.7,-0.2){\scriptsize $6$} \put(-0.2,1.5){\scriptsize $7$}
\end{picture}
$\subseteq {\rm Lk}_X(25)$, a contradiction to Observation 2.

\medskip

\noindent {\sf Case 11\,:} $U=\Sigma_{\,4}$ of Example \ref{e1}.
Thus, the odd triangles are $124$, $127$, $145$, $156$, $167$,
$234$, $237$, $347$, $457$ and $567$. By the claim, $1247$,
$1457$, $1567$, $2347 \in X$. Then the sum of the remaining
sixteen tetrahedra gives a non-zero element of $H_3(X, \ZZ_2)$, a
contradiction.

\medskip

\noindent {\sf Case 12\,:} $U=\Sigma_{\,5}$ of Example \ref{e1}.
Thus, the odd triangles are  $123$, $126$, $135$, $156$, $234$,
$246$, $345$, $457$, $467$, $567$. By the claim, $1234$, $1235$,
$1246$, $1256$, $1345$, $2345$, $3457$, $4567$ $\in X$. Thus ${\rm
Lk}_X(14)\supseteq $ \setlength{\unitlength}{2mm}
\begin{picture}(4.5,2.7)(0,-0.1)

\put(0.7,0){$_{\bullet}$} \put(0.7,2){$_{\bullet}$}
\put(2.7,0){$_{\bullet}$} \put(2.7,2){$_{\bullet}$}

\thinlines

\put(3,0){\line(0,1){2}} \put(1,0){\line(1,0){2}}
\put(1,2){\line(1,0){2}}

\put(3.7,1.5){\scriptsize $2$} \put(-0.2,-0.2){\scriptsize $5$}
\put(3.7,-0.2){\scriptsize $3$} \put(-0.2,1.5){\scriptsize $6$}
\end{picture} and ${\rm Lk}_X(25)\supseteq $
\setlength{\unitlength}{2mm}
\begin{picture}(4.5,2.7)(0,-0.1)

\put(0.7,0){$_{\bullet}$} \put(0.7,2){$_{\bullet}$}
\put(2.7,0){$_{\bullet}$} \put(2.7,2){$_{\bullet}$}

\thinlines

\put(3,0){\line(0,1){2}} \put(1,0){\line(1,0){2}}
\put(1,2){\line(1,0){2}}

\put(3.7,1.5){\scriptsize $1$} \put(-0.2,-0.2){\scriptsize $4$}
\put(3.7,-0.2){\scriptsize $3$} \put(-0.2,1.5){\scriptsize $6$}
\end{picture}.
Since $14$ and $25$ are not in $U$, Observation 2 implies that
${\rm Lk}_X(14) = $ \setlength{\unitlength}{2mm}
\begin{picture}(6.6,3)(0,0)

\put(0.9,1){$_{\bullet}$} \put(2.7,0){$_{\bullet}$}
\put(2.7,2){$_{\bullet}$} \put(4.7,0){$_{\bullet}$}
\put(4.7,2){$_{\bullet}$}

\thinlines

\put(3,0){\line(1,0){2}} \put(5,0){\line(0,1){2}}
\put(3,2){\line(1,0){2}} \put(1,1){\line(2,1){2}}
\put(1,1){\line(2,-1){2}}

\put(1.8,2){\scriptsize 6} \put(0,-0.2){\scriptsize 7}
\put(3.4,0.4){\scriptsize 5} \put(5.7,-0.1){\scriptsize 3}
\put(5.7,2){\scriptsize 2}
\end{picture}
and ${\rm Lk}_X(25) = $ \setlength{\unitlength}{2mm}
\begin{picture}(6.6,3)(0,0)

\put(0.9,1){$_{\bullet}$} \put(2.7,0){$_{\bullet}$}
\put(2.7,2){$_{\bullet}$} \put(4.7,0){$_{\bullet}$}
\put(4.7,2){$_{\bullet}$}

\thinlines

\put(3,0){\line(1,0){2}} \put(5,0){\line(0,1){2}}
\put(3,2){\line(1,0){2}} \put(1,1){\line(2,1){2}}
\put(1,1){\line(2,-1){2}}

\put(1.8,2){\scriptsize 6} \put(0,-0.2){\scriptsize 7}
\put(3.4,0.4){\scriptsize 4} \put(5.7,-0.1){\scriptsize 3}
\put(5.7,2){\scriptsize 1}
\end{picture}.
Thus $1457$, $2457\in X$. Then the triangle $457$ is of degree 4
in $X$, a contradiction. This completes the proof. \hfill $\Box$

\bigskip

\noindent {\bf Proof of Theorem \ref{t1}\,.} Let $Y$ be a minimal
counter example. So, $Y$ is an $n$-vertex (for some $n\leq 7$)
$\ZZ_2$-acyclic simplicial complex which is not collapsible to any
proper subcomplex.

If $n < 7$ then choose a facet $\alpha$ of $Y$ and an element
$v\not\in V(Y)$. Let $\widetilde{Y}$ be obtained from $Y$ by the
bistellar $d$-move $\kappa_{\alpha \cup \{v\}}$, where $d$ is the
dimension of $Y$. Then $\widetilde{Y}$ is an $(n+1)$-vertex
$\ZZ_2$-acyclic simplicial complex. Since $Y$ has no free face,
$\widetilde{Y}$ has no free face and hence $\widetilde{Y}$ is not
collapsible to any proper subcomplex. Repeating this construction
(if necessary) we get a 7-vertex $\ZZ_2$-acyclic simplicial
complex $X$ which is not collapsible to any proper subcomplex.
Then, by Lemma \ref{l3.1}, $X$ is of dimension $2$ or $3$. But,
this is not possible by Lemmas \ref{l3.2} and \ref{l3.5}. This
completes the proof. \hfill $\Box$


\section{Homology spheres.}

A connected $d$-dimensional weak pseudomanifold is called a {\em
normal pseudomanifold} if the links of all the simplices of
dimension up to $d - 2$ are connected. Observe that if $X$ is
a normal pseudomanifold then $X$ is a pseudomanifold. (If not
then, since $X$ is connected, there exist two intersecting
facets $\tau$, $\sigma$ for which there is no sequence of
facets $\tau = \tau_0, \dots, \tau_n=\sigma$ such that $\tau_{i-1}
\cap \tau_{i}$ is a face of codimension 1 for $1\leq i\leq n$.
Choose $\tau$, $\sigma$ among all
such pairs such that $\dim(\tau \cap \sigma)$ is maximum.
Then $\dim(\tau \cap \sigma) \leq d-2$ and ${\rm
lk}_X(\tau \cap \sigma)$ is not connected, a contradiction.)
Notice that all the links of positive dimensions (i.e., the links
of simplices of dimension $\leq d - 2$) in a normal
$d$-pseudomanifold are normal pseudomanifolds (and hence are
pseudomanifolds). Clearly, any triangulation of a
connected closed manifold is a normal pseudomanifold.

\begin{lemma}$\!\!${\bf .} \label{l4.1}
Let $Y$ be a $d$-dimensional normal pseudomanifold. Let $Y_1$ be a
proper induced subcomplex of $Y$ which is pure of dimension $d$.
Put $L=C(Y_1, Y)$ and $Y_2= N(L, Y)$. Then $(a)$ $Y_1$, $Y_2$ are
weak pseudomanifolds with boundary, $(b)$ $\partial Y_2$ is an
induced subcomplex of $Y_2$ and $(c)$ $\partial Y_2 = \partial Y_1
= Y_1\cap Y_2$.
\end{lemma}

\noindent {\bf Proof\,.} Since $Y$ is a pseudomanifold and
$Y_1\subset Y$ is pure of maximum dimension, $Y_1$ is a weak
pseudomanifold with boundary.  Since the maximal simplices of
$Y_2$ are those maximal simplices of $Y$ which intersect $V(L)$,
$Y_2$ is pure of dimension $d$ and each $d$-simplex of $Y$ is
either in $Y_1$ or in $Y_2$ but not in both. This implies that
$Y_2$ is a weak pseudomanifold with boundary. This proves $(a)$.

Let $V_1= V(Y_1)$, $V_2=V(L)$. Then $V(Y) = V_1 \sqcup V_2$. Now,
$\tau$ is a facet of $\partial Y_2$ $\Leftrightarrow$ there exists
a unique $d$-face $\sigma_2\in Y_2$ containing $\tau$
$\Leftrightarrow$ there exists a unique $d$-face $\sigma_1\in Y_1$
containing $\tau$ $\Leftrightarrow$ $\tau$ is a facet of $\partial
Y_1$. Therefore, $\partial Y_2 = \partial Y_1\subseteq Y_1\cap
Y_2$.

Conversely, let $\tau$ be an $i$-simplex in $Y_1 \cap Y_2$. If
possible, let $\tau \not\in \partial Y_1 = \partial Y_2$. Then
$\tau \in Y_1 \setminus \partial Y_1$. Therefore, ${\rm
lk}_{Y_1}(\tau)$ is a $(d - i - 1)$-dimensional weak
pseudomanifold (without boundary) and ${\rm lk}_{Y_1}(\tau)
\subseteq {\rm lk}_Y(\tau)$. If $i = d - 1$ then each of ${\rm
lk}_{Y_1}(\tau)$ and ${\rm lk}_{Y}(\tau)$ consists of two vertices
and hence ${\rm lk}_{Y_1}(\tau) = {\rm lk}_{Y_1}(\tau)$. Now,
assume that $i \leq d - 2$. Since $Y$ is a normal pseudomanifold,
${\rm lk}_Y(\tau)$ is a pseudomanifold (of dimension $\geq 1$) and
hence ${\rm lk}_{Y_1}(\tau) = {\rm lk}_Y(\tau)$. This implies that
${\rm star}_Y(\tau) \subseteq Y_1$ and hence ${\rm star}_Y(\tau)
\cap L = \emptyset$. Now, $\tau \in Y_2$ and $Y_2$ is pure.
Therefore, there exists a $d$-simplex $\sigma\in Y_2$ such that
$\tau \subseteq \sigma$. Since $\sigma \in Y_2$, $\sigma \cap V(L)
= \emptyset$. This implies that ${\rm star}_Y(\tau) \cap L \neq
\emptyset$, a contradiction. Therefore, $\tau \in \partial Y_1 =
\partial Y_2$. So, $Y_1\cap Y_2 = \partial Y_1 = \partial Y_2$.
This proves (c).

Since $\partial Y_2 = \partial Y_1$, $\partial Y_2\subseteq
Y_2[V_1] = Y_2[V_1]\cap Y[V_1]= Y_2[V_1]\cap Y_1\subseteq Y_2\cap
Y_1 = \partial Y_2$. Thus, $\partial Y_2 = Y_2[V_1] = Y_2[V_1\cap
V(Y_2)]$. This proves $(b)$. \hfill $\Box$

\begin{lemma}$\!\!${\bf .} \label{l4.2}
Let $X$ be a connected combinatorial $d$-manifold. Let $X_1$ be an
induced subcomplex of $X$ which is a combinatorial $d$-ball. Put
$L = C(X_1, X)$ and $X_2= N(L, X)$. Then
\begin{enumerate}
     \item[$(a)$] $X_2$ is a connected combinatorial $d$-manifold with
     boundary.
     \item[$(b)$] $|X_2| \coll |L|$.
     \item[$(c)$] If, further, $L$ is collapsible then $X$ is a
     combinatorial sphere.
     \end{enumerate}
\end{lemma}

\noindent {\bf Proof\,.} Let $V_1= V(X_1)$, $V_2=V(L)$. Then $V(X)
= V_1 \sqcup V_2$. As in the proof of Lemma \ref{l4.1}, $X_2$ is
pure of dimension $d$ and each $d$-simplex of $X$ is either in
$X_1$ or in $X_2$ but not in both.

Let $v$ be a vertex of $X_2$. Notice that $v \in X_1\setminus
\partial X_1$ $\Rightarrow$ ${\rm Lk}_{X_1}(v)\subseteq {\rm
Lk}_{X}(v)$ are $(d - 1)$-spheres $\Rightarrow$ ${\rm Lk}_{X_1}(v)
= {\rm Lk}_{X}(v)$ $\Rightarrow$ $v\not\in X_2$, a contradiction.
So, either $v\in V_2$ or $v\in \partial X_1$.

If $v \in V_2$ then each $d$-simplex of $X$ containing $v$ is in
$X_2$ and hence ${\rm Lk}_{X_2}(v)= {\rm Lk}_{X}(v)$ is a
combinatorial $(d-1)$-sphere.

If $v\in \partial X_1$ then $(Y, Y_1, Y_2) := ({\rm Lk}_X(v), {\rm
Lk}_{X_1}(v), {\rm Lk}_{X_2}(v))$ satisfies the hypothesis of
Lemma \ref{l4.1}. Therefore, by Lemma \ref{l4.1}, ${\rm
Lk}_{X_1}(v) \cap {\rm Lk}_{X_2}(v) = \partial ({\rm
Lk}_{X_2}(v))$. This implies that the closure of $|{\rm
Lk}_{X}(v)| \setminus |{\rm Lk}_{X_1}(v)|$ in $|{\rm Lk}_{X}(v)|$
is $|{\rm Lk}_{X_2}(v)|$. Since $|{\rm Lk}_{X}(v)|$ is a pl
$(d-1)$-sphere and  $|{\rm Lk}_{X_1}(v)|$ is a pl $(d-1)$-ball, by
Proposition \ref{ks3.13}, $|{\rm Lk}_{X_2}(v)|$ is a pl $(d
-1)$-ball. Thus, ${\rm Lk}_{X_2}(v)$ is a combinatorial
$(d-1)$-ball.

Thus $X_2$ is a combinatorial $d$-manifold with boundary such that
$\partial X_2$ ($= \partial X_1$, by Lemma \ref{l4.1}) is
connected. Therefore, if $X_2$ were disconnected, it would have a
$d$-dimensional weak pseudomanifold as a component. This is not
possible since $X$ is a $d$-dimensional pseudomanifold. Therefore
$X_2$ is connected. This proves $(a)$.

As $L = X[V_2]$, we have $L\subseteq X_2$ and hence $L= X_2[V_2]$.
Since, by Lemma \ref{l4.1}, $\partial X_2$ is the induced
subcomplex of $X_2$ on $V_1\cap V(X_2)$, this implies that $L$ is
the simplicial complement of $\partial X_2$ in $X_2$. Then, by
Proposition \ref{snt}, $|X_2| \coll |L|$. This proves $(b)$.

Now, if $L \scoll 0$ then $|L|\coll 0$ and hence $|X_2|\coll 0$.
So, by Proposition \ref{ks3.28}, $|X_2|$ is a   pl ball.

Let $\sigma$ be a $d$-simplex in $S^{\,d}_{d + 2}$. Let $B_1 =
|\sigma|$ and $B_2 = |S^{\,d}_{d + 2}\setminus\{\sigma\}|$. Then
$B_1$ and $B_2$ are pl $d$-balls. Let $f_2 \colon B_2 \to |X_2|$
be a pl homeomorphism. Let $f = f_2|_{\partial B_2}$. Since
$\partial B_1 = \partial B_2$ and $\partial (|X_1|) = |\partial
X_1| = |\partial X_2|$, $f \colon \partial B_1 \to \partial
(|X_1|)$ is a pl homeomorphism. By Proposition \ref{ks1.10}, there
exists a pl homeomorphism $f_1 \colon B_1 \to |X_1|$ such that
$f_1|_{\partial B_1} = f = f_2|_{\partial  B_2}$. Then $f_1\cup
f_2$ is a pl homeomorphism from $|S^{\,d}_{d+2}|$ to $|X|$. This
proves $(c)$. \hfill $\Box$

\begin{lemma}$\!\!${\bf .} \label{l4.3}
Let $X$ be a combinatorial triangulation of a $\ZZ_2$-homology
$d$-sphere. Let $X_1$ be an induced subcomplex of $X$ which is a
combinatorial $d$-ball. Let $L = C(X_1, X)$ and $X_2 = N(L, X)$.
Then $X_2$ is $\ZZ_2$-acyclic.
\end{lemma}

\noindent {\bf Proof\,.} Let $J=X_1\cap X_2$. Then, by Lemma
\ref{l4.1}, $J = \partial X_1$. So, $J$ is a combinatorial
$(d-1)$-sphere. Therefore, $H_{d-1}(J, \ZZ_2) = \ZZ_2$ and
$\widetilde{H}_q(J, \ZZ_2)=0$ for all $q \neq d-1$. Also
$\widetilde{H}_q(X_1, \ZZ_2)=0$ for all $q\geq 0$. For $q\geq 1$,
we have the following exact Mayer-Vietoris sequence of homology
groups with coefficients in $\ZZ_2$ (see \cite{m, sp})\,:
\begin{equation} \label{mvs}
\cdots \to H_{q+1}(X) \to H_{q}(J) \to H_{q}(X_1) \oplus
H_{q}(X_2) \to H_{q}(X)\to \widetilde{H}_{q-1}(J) \to \cdots
\end{equation}

Now, $H_d(X, \ZZ_2)=\ZZ_2$ and $\widetilde{H}_q(X, \ZZ_2)= 0$ for
$q \neq d$. By Lemma \ref{l4.2}, $|X_2|$ is a connected
$d$-manifold with non-trivial boundary. Therefore, $H_d(X_2,
\ZZ_2)=0$ and $H_0(X_2, \ZZ_2) = \ZZ_2$. Then, by (\ref{mvs}),
$H_q(X_2, \ZZ_2) = 0$ for $0 < q < d - 1$ and for $q = d - 1$ we
get the following short exact sequence of abelian groups\,:
$$
0 \to \ZZ_2 \to \ZZ_2 \to H_{d - 1}(X_2, \ZZ_2) \to 0.
$$
Clearly, this implies $H_{d - 1}(X_2, \ZZ_2)=0$.  Thus,
$\widetilde{H}_{q}(X_2, \ZZ_2)=0$ for all $q\geq 0$. \hfill $\Box$

\bigskip

\noindent {\bf Proof of Theorem \ref{t2}.} Let $X_1$ be an
$m$-vertex induced subcomplex of $M$ which is a combinatorial
$d$-ball. Let $L = C(X_1, M)$ and $X_2= N(L, M)$. Then, by Part
$(b)$ of Lemma \ref{l4.2}, $|X_2|\coll |L|$.

Again, by Lemma \ref{l4.3}, $X_2$ is $\ZZ_2$-acyclic and hence $L$
is $\ZZ_2$-acyclic. Since $n\leq m+7$, the number of vertices in
$L$ is $\leq 7$. Therefore, by Theorem \ref{t1}, $L$ is
collapsible. Then, by Part $(c)$ of Lemma \ref{l4.2}, $M$ is a
combinatorial sphere. \hfill $\Box$

\bigskip

\noindent {\bf Proof of Corollary \ref{t3}.} If $\sigma$ is a
$d$-simplex of $M$ then the induced subcomplex
$\Delta^d_{d+1}(\sigma)$ is a $(d+1)$-vertex combinatorial
$d$-ball. Therefore, by Theorem \ref{t2}, $M$ is a combinatorial
sphere. \hfill $\Box$

\medskip

\noindent {\bf Proof of Corollary \ref{t4}.} Assume, if possible,
that $M$ admits a bistellar $i$-move $\kappa_A$ for some $i < d$.
Let $\beta$ be the core of $A$ and $\alpha = A\setminus \beta$.
Then $M[A] = \Delta^i_{i+1}(\alpha)\ast S^{\,d-i-1}_{d-i+
1}(\beta)$ is a $(d+2)$-vertex combinatorial $d$-ball. Therefore,
by Theorem \ref{t2}, $M$ is a combinatorial sphere, a
contradiction. This proves the corollary. \hfill $\Box$

\medskip


{\footnotesize

\begin{thebibliography}{99}
\bibitem{bd2}
B. Bagchi, B. Datta, A structure theorem for pseudomanifolds, {\em
Discrete Math.} {\bf 168} (1998) 41--60.
\bibitem{bd4}
B. Bagchi, B. Datta, Non-existence of 6-dimensional
pseudomanifolds with complementarity, {\em Adv. Geom.} {\bf 4}
(2004), 537--550.
\bibitem{b}
R. H. Bing, {\em The Geometric Topology of $3$-Manifolds}, Amer.
Math. Soc., Providence, 1983.
\bibitem{bl}
A. Bj\"{o}rner, F. H. Lutz,  Simplicial manifolds, bistellar flips
and a 16-vertex triangulation of the Poincar\'{e} homology
3-sphere, {\em Exp. Math} {\bf 9} (2000) 275--289.
\bibitem{bk}
U. Brehm, W. K\"{u}hnel, Combinatorial manifolds with few
vertices, {\em Topology} {\bf 26} (1987) 465--473.
\bibitem{d}
B. Datta, Two dimensional weak pseudomanifolds on seven vertices,
{\em Bol. Soc. Mat. Mexicana} {\bf 5} (1999) 419--426.
\bibitem{k}
W. K\"{u}hnel, {\em Tight Polyhedral Submanifolds and Tight
Triangulations}, Lecture Notes in Mathmatics, {\bf 1612},
Springer-Verlag, Berlin, 1995.
\bibitem{l}
F. H. Lutz, {\em Triangulated Manifolds with Few Vertices and
Vertex-Transitive Group Actions}, Thesis (TU, Berlin), Shaker
Verlag, Aachen, 1999.
\bibitem{m}
J. R. Munkres, {\em Elements of Algebraic Topology\/},
Addison-Wesley, Menlo Park, CA, 1984.
\bibitem{p}
U. Pachner, Konstruktionsmethoden und das kombinatorische
Hom\"{o}omorphieproblem f\"{u}r Triangulationen kompakter
semilinearer Mannigfaltigkeiten, {\em Abh. Math. Sem. Univ.
Hamburg} {\bf 57} (1987) 69--86.
\bibitem{rs}
C. P. Rourke, B. J. Sanderson, {\em Introduction to
Piecewise-Linear Topology}, Springer-Verlag, Berlin, 1982.
\bibitem{sp}
E. H. Spanier, {\em Algebraic Topology}, Springer-Verlag, New
York, 1966.
\bibitem{st}
R. P. Stanley, A combinatorial decomposition of acyclic simplicial
complexes, Discrete Math. {\bf 120} (1993) 175--182.
\bibitem{w}
D. W. Walkup, The lower bound conjecture for 3- and 4-manifolds,
{\em Acta Math.} {\bf 125} (1970) 75--107.
\end{thebibliography}
}
\end{document}